\documentclass[12pt,reqno]{amsart}
\usepackage{amsmath}
\usepackage{amssymb}

\usepackage{psfrag}
\usepackage{graphicx}

\begingroup
\newtheorem{theorem}{Theorem}[section]

\newtheorem{proposition}[theorem]{Proposition}

\endgroup

\theoremstyle{definition}
\begingroup
\newtheorem{definition}[theorem]{Definition}
\newtheorem{remark}[theorem]{Remark}

\endgroup
\numberwithin{equation}{section}

\newcommand{\bbibitem}{\bibitem}

\newcommand{\llabel}[1]{\label{#1}}

\newcommand{\R}{{\mathbb R}}
\newcommand{\e}{{\varepsilon}}
\newcommand{\E}{{\mathcal{E}}}
\newcommand{\W}{{\mathcal{W}}}
\newcommand{\F}{{\mathcal{F}}}
\newcommand{\GG}{{\mathcal{G}}}

\newcommand{\G}{{\Gamma}}
\newcommand{\hh}{{\mathcal{H}^1}}
\newcommand{\kk}{{\kappa}}

\newcommand{\sg}{{\sigma}}
\newcommand{\Om}{{\Omega}}
\newcommand{\Omo}{{\Omega\setminus\Gamma(\sg_0)}}

\newcommand{\grad}{{\rm grad}}
\newcommand{\wto}{{\rightharpoonup\,}}
\newcommand{\ue}{{u^\e}}

\title[A model of quasistatic crack growth]{An artificial viscosity approach \\[3mm]
to quasistatic crack growth }

\author[Rodica Toader]{Rodica Toader}\address[Rodica Toader]{Universit\`a degli Studi di Udine, Dipartimento di 
Ingegneria Civile, Via delle Scienze 208, 
33100 Udine, Italy}\email[Rodica Toader]{toader@uniud.it}

\author[Chiara Zanini]{Chiara Zanini}
\address[Chiara Zanini]{SISSA, Via Beirut 4, 34014 Trieste, Italy}\email[Chiara Zanini]
{zaninic@sissa.it}

\begin{document}

\footnotetext{Preprint SISSA 43/2006/M (July 2006).}

\noindent
\begin{abstract}
\noindent We introduce a new model of irreversible quasistatic crack growth in which the evolution of cracks is the limit of a suitably modified $\e$-gradient flow of the energy functional, as the  ``viscosity" parameter $\e$ tends to zero. 
\end{abstract}

\maketitle

 {\small\bigskip
}

\begin{section}{Introduction} \llabel{introd}

In this paper we consider the crack growth in brittle materials in the particular case 
of a preassigned crack path $\Gamma$. We assume that $\Gamma$ is a regular arc with one endpoint on the boundary of the reference configuration $\Om$, %and the other in the interior of $\Om$. Moreover, we assume 
that there exists an initial connected crack starting from the boundary point, and that the crack remains connected during the evolution. Hence, such a  crack will be completely determined by its length $\sg$.
The evolution is assumed to be irreversible, so that the length of the crack will be increasing in time,
and  quasistatic, i.e. at each time the configuration describing the body is in equilibrium.    By {\em configuration} we mean a pair $(u,\sg)$ where $u$ represents the displacement, and $\sg$ the length of the crack.

The main new feature of this model is that it is based on a local stability criterion for the energy functional rather than on a global one. The choice of the total energy $\E(t)(u,\sg)$, of a configuration $(u,\sg)$ at time $t$, is inspired by Griffith's idea \cite{GR} that the evolution of cracks in brittle materials is the result of the competition between the elastic energy of the body and the energy needed to extend the crack. For simplicity of exposition, we detail our study in  the case of antiplane shears and therefore, $u$ is a scalar function representing the displacement orthogonal to the plane of $\Om$ and the bulk part of the energy is given by the square of the $L^2$-norm of the gradient of $u$. We consider the case of a homogeneous isotropic material and, according to Griffith's theory \cite{GR}, we assume the surface energy to be proportional to the length $\sg$ of the crack, the constant of proportionality being given by the toughness of the material. The evolution is driven by time-dependent imposed boundary displacements $\psi(t)$ on a part $\partial_D\Om$ of the boundary, and applied boundary forces $g(t)$ on the remaining part $\partial_N\Om$. The total energy, $\E(t)(u,\sg)$, is the sum of the bulk energy and the surface energy minus the work of the applied forces $g(t)$. 

Let $AD(\psi(t),\sg)$ denote the set of {\em admissible displacements},  i.e.  displacements with finite bulk energy, compatible with the imposed boundary displacement $\psi(t)$ and with the crack length $\sg$. Let us recall that the functional $\E(t)(u,\sg)$ is not differentiable, nor convex. It depends on $\sg$ both through the surface energy term and through the constraint on the set of admissible displacements. Note that, given $t$ and $\sg$, there exists a unique  minimizer $u_{t,\sg}$ of the energy $\E(t)(u,\sg)$ in $AD(\psi(t),\sg)$. Then let us consider the  minimal energy corresponding to the boundary data $\psi(t)$ and to the crack length $\sg$: $E(t,\sg):=\E(t)(u_{t,\sg},\sg)$.  
The derivative $\partial_\sg E(t,\sg)$ can be computed (see Proposition \ref{p:ensg}) and it is related to the stress intensity factor of the displacement $u_{t,\sg}$ at the tip of the crack. It plays a crucial r\^ole in the Griffith's criterion for the propagation of cracks.

Let us define now the notion of evolution we are interested in.  
The {\em irreversible quasistatic evolution problem} 
consists in finding  a
left-continuous function of time
$t\mapsto(u(t),\sg(t))$  
such that the displacement $u(t)$ at time $t$  belongs to the set $AD(\psi(t),\sg(t))$, and  the following
three conditions are satisfied:  
\begin{itemize}
\item[(a)] {\em local unilateral stability:} at every time $t\geq0$
\begin{eqnarray*}
&& \E(t)(u(t),\sg(t))\leq \E(t)(v,\sg(t)) \qquad\forall\, v\in AD(\psi(t),\sg(t)) \\
&&\partial_\sg E(t,\sg(t))\geq 0\,,
\end{eqnarray*}
\item[(b)] {\em irreversibility:} the map $t\mapsto \sg(t)$ is increasing;
\item[(c)] {\em energy inequality:} for every $0\leq s<t$ we have
\begin{equation*}
\E(t)(u(t),\sg(t)) \leq \E(s)(u(s),\sg(s))+ {\rm Work}(u;s,t)\,,
\end{equation*}
where ${\rm Work}(u;s,t)$ denotes the work of external forces.
\end{itemize} 
A solution to this problem will be called an {\em irreversible quasistatic evolution}.

We show that conditions (a)-(c) are enough to ensure that at almost every time  $t$ a weak version of Griffith's criterion is satisfied (see Proposition 2.6). In condition (c) two terms contribute to the work of the external forces: the first is due to the surface forces generated by the imposed boundary displacement and the second one comes from the applied surface loads.

In \cite{FrMa-98}, \cite{DM-T-02}, \cite{C}, \cite{F-L}, \cite{DM-F-T} the  {\em globally stable irreversible quasistatic evolution problem} was studied. It is a particular case of the previous one and it fits the
general scheme of the continuous-time energetic formulation of rate-independent
processes developed by Mielke and his collaborators (see \cite{M} and the references therein).
It  consists in  
finding an irreversible quasistatic evolution
which satisfies the {\em global stability condition:} at every time $t\geq0$
\begin{equation*}
\E(t)(u(t),\sg(t))\leq \E(t)(v,\sg)\qquad\forall\,\sg\geq\sg(t)\quad\forall\, v\in
AD(\psi(t),\sg)\,.
\end{equation*}
In this case condition (c) can be replaced by the {\em energy balance} condition: the increment in stored energy plus the energy spent in crack increase equals the work of external forces.  

The  global minimality condition imposes the comparison, in terms of energy, of a configuration with all admissible configurations with a longer crack  and might generate jumps in the length of the crack that are not justified by the mechanical interpretation of the problem. That is why we look for a selection criterion different from the global stability. 

The selection criterion we choose is based on an approximation procedure with a regularizing effect. We consider solutions $(u_\e(t),\sg_\e(t))$ of regular evolution problems that converge to a solution $(u(t),\sg(t))$ of the irreversible quasistatic evolution problem. Moreover, among the possible approximations we choose one that has the following property:  
\begin{itemize}
\item[$({\mathcal P})$]
if on a certain time interval $[t_1,t_2]$ there exists a regular function $\sg_0(t)$ such that
\begin{equation*}
\partial_\sg E(t,\sg_0(t))=0\quad\hbox{and}\quad
\partial^2_\sg E(t,\sg_0(t))>0\qquad\forall\,t\in[t_1,t_2]\,,
\end{equation*}
and if $\dot\sg_\e(t)>0$ for every $t\in[t_1,t_2]$, then the equality $\sg(t_1)=\sg_0(t_1)$ implies that $\sg(t)=\sg_0(t)$ for every $t\in[t_1,t_2]$.
\end{itemize}

Let us now describe more in detail the construction of the approximating evolutions. First of all, we fix an initial condition: assume that at time $t=0$ the crack length is equal to 
$\sg_0>0$ and  the displacement is equal to 
$u_0$, in such a way that the initial configuration $(u_0,\sg_0)$ is in equilibrium. 
Then, for every $\sg$ between $\sg_0$ and $\overline\sg$, where $\overline\sg$ is the length of $\Gamma$,  we consider a  diffeomorphism $\Phi_\sg$ of $\Om$ that transforms the crack of length $\sg$ into the one of length $\sg_0$.
Using $\Phi_\sg$, we change 
variables in the expression of the energy functional
$\E$ and transform it into a functional
$\F$ depending on the time $t$, the crack length $\sg$, and  the modified displacement $v$, which, assuming that the change of variables does not influence the terms due to the boundary data, has the form
\begin{eqnarray*}
\F(t,v,\sg)=\int_\Omo(A(\sg,x)(Dv(x)+D\psi(t,x))|Dv(x)+D\psi(t,x))dx+\sg-\\
-\int_{\partial_N\Om} g(t,x)(v(x)+\psi(t,x))d\hh(x)\,.
\end{eqnarray*}
Here, $\Gamma(\sg_0)$ is the crack of length $\sg_0$, $A(\sg,x)$ is a $2{\times}2$ symmetric matrix of smooth coefficients coming from the change of variables, $Dv$ is the distributional gradient of $v$ with respect to the spatial variables $x$, $(\cdot|\cdot)$ denotes the scalar product in $\R^2$, and $\hh$  the one-dimensional Hausdorff measure. 
The advantage of this change of variables is that now  the set of admissible functions $v$ does not depend on $t$, nor on the crack length $\sg$. The same change of variables is considered, in a suitable small neighbourhood of the crack tip, in order to compute the derivative $\partial_\sg E(t,\sg)$ (see also \cite{Gr-92}, \cite{CD-DM-00}, \cite{K-M}).

The approximating evolution $(v_\e,\sg_\e)$ is the solution of a suitably
modified
$\e$-gradient flow for the functional $\F$ which starts from the initial data
$(u_0,\sg_0)$:
\begin{equation}\label{idefegradflow}
\begin{cases}
\e\dot v_\e=-\grad_v \F(t,v_\e,\sg_\e)\\
\displaystyle\e\dot\sg_\e=(-\partial_\sg \F(t,v_\e,\sg_\e))^+\lambda(\sg_\e)\\ 
v_\e(0)=u_0\\
\sg_\e(0)=\sg_0\,.
\end{cases}
\end{equation}
Here $\grad_v\F(t,v,\sg)$ denotes the gradient of the function $v\mapsto \F(t,v,\sg)$ considered as a function defined on the Sobolev space $H^1(\Omo)$ with suitable boundary conditions. The positive part in the second equation guarantees the irreversibility of the evolution, while $\lambda$ is a Lipschitz continuous positive cut-off function that becomes zero for $\sg=\overline\sg$, so that only increasing solutions with crack length less than $\overline\sg$ are considered. 
If we are interested in the evolution %of the configurations $(u,\sg)$ 
until a certain crack length $\sg_1$, with $\sg_0<\sg_1<\overline\sg$, is reached,  then we choose $\lambda(\sg)=1$ for $\sg_0\leq\sg\leq\sg_1$. In this way for crack lengths less than $\sg_1$, the regularized evolution law  is proportional to the gradient flow for $\F$, with the constraint that the crack length is increasing, while it is 
 distorted by $\lambda$
 for crack lengths between $\sg_1$ and  $\overline\sg$. 
 Therefore the evolution is considered meaningful only until the crack reaches the length $\sg_1$. 

Note that, using the form of the functional $\F$, the first equation in (\ref{idefegradflow}) can be written as
$$
\e\Delta_x\dot v_\e(t,x)=-{\rm div}_x(A(\sg_\e(t),x)Dv_\e(t,x))+\dots
$$
with suitable boundary conditions.
We preferred the evolution problem in $H^1$ to the usual parabolic one
$$
\e\dot v_\e(t,x)=-{\rm div}_x(A(\sg_\e(t),x)Dv_\e(t,x))+\dots
$$
which corresponds to the gradient flow in $L^2$, because it helped us to prove property~$({\mathcal P})$, see Theorem~\ref{tfi}. 
Note  also that in this way the first equation in (\ref{idefegradflow}) becomes an ODE and thus the existence of the solution for this modified 
$\e$-gradient flow follows from classical existence and uniqueness results for ordinary differential equations in Banach spaces.

Let us remark that this model is not suited for the study of the crack initiation problem.   We also note that the approximating evolutions we consider have been chosen on the basis of their mathematical simplicity and do not seem to have any  mechanical interpretation. Nevertheless, we think that the notion of {\em approximable irreversible quasistatic evolution} proposed here could be the starting point for the study of different approximations with a mechanical justification.
For a different approach to the irreversible quasistatic crack growth see also~\cite{FHV}.

 In Section 5 we detail our results in the case of monotonically increasing in time imposed boundary displacements and compare this evolution with the one proposed by Francfort and Marigo in \cite{FrMa-98}, while in Section 6 we provide an example where the energy, as function of the crack length, has at least a concavity interval.

\end{section}

\begin{section}{Setting of the problem}

\begin{subsection}{The reference configuration and the crack.}  
  Let $\Om$ be a bounded connected open set of $\R ^2$ with Lipschitz boundary 
$\partial\Omega$. The set $\overline\Om$ represents the reference configuration
of an isotropic, homogeneous elastic body. Let 
$\partial_D\Om$ be a closed subset of $\partial\Om$ with $\hh(\partial_D\Om)>0$, where $\hh$ denotes the one-dimensional Hausdorff measure, and let 
$\partial_N\Om:=\partial\Om\setminus\partial_D\Om$. On the Dirichlet part of the
boundary, $\partial_D\Om$, we will impose the boundary displacements, while on the
Neumann part of the boundary, $\partial_N\Om$, we will prescribe the boundary forces.

Let $\Gamma$ be a simple $C^3$-arc and let $\gamma\colon [0,\overline \sg]\to \G$ be its arc-length parametrization. We assume that  $\gamma(0)\in\partial_N\Om$ and  $\gamma( \sg)\in \Om$ for $0<\sg\leq\overline\sg$. For technical reasons it is convenient to extend $\Gamma$ until it reaches another point in $\partial_N\Om$, so that it cuts the reference configuration $\Omega$ into two subsets. The extension will still be called $\Gamma$, and its arc-length parametrization will now be $\gamma\colon [0,\sg_{max}]\to \G$. We assume that its intersection with the boundary $\partial\Omega$ is not tangential. Let $\nu$ be a unit normal  vector field on $\Gamma$. Then we denote by $\Omega^+$ the part of $\Omega\setminus \G$ which is positively oriented with respect to $\nu$, and by $\Omega^-$ the remaining part, so that
$\Omega\setminus\G=\Omega^+\cup\Omega^-$. Both $\Omega^+$ and $\Omega^-$
are bounded connected sets with Lipschitz boundary. We assume that 
$\hh(\partial_D\Omega\cap \partial\Omega ^+)>0$ and  
$\hh(\partial_D\Omega\cap \partial\Omega ^-)>0$.
We make the following
simplifying assumption: all admissible cracks are of the form
\begin{equation*}
\G(\sg):= \{\gamma(s): 0\leq s\leq \sg\}\qquad\hbox{with }\sg\leq\overline \sg\,.
\end{equation*}

According to Griffith's theory we assume that the energy spent to produce the
crack $\G(\sg)$ is proportional to the length of the crack, and, for simplicity,
we take it to be equal to $\sg$.

\end{subsection}

\begin{subsection}{ The bulk energy.} We consider here the case of antiplane shears.
Given a crack $\G(\sg)$, an admissible displacement is any function $u\in
H^1(\Om\setminus\G(\sg))$, and the bulk energy associated to the
displacement
$u$ is 
$$
\W(D u):=\int_{\Om\setminus\G(\sg)}|D u(x)|^2dx\,,
$$
where $D u$ is the distributional gradient of $u$ and $|\cdot|$ denotes the norm in $\R ^2$.
\end{subsection}

\begin{subsection}{ The boundary displacement } In the following it will be convenient to work on a fixed time interval $[0,T]$ with $T>0$. We impose a time-dependent
Dirichlet boundary condition on $\partial_D\Om$: 
\begin{equation*}
u=\psi(t) \quad \text{\ on $\partial_D\Omega$},
\end{equation*}
where equality on the boundary is considered in the sense of traces.
We assume that $\psi(t)$ is the trace on $\partial_D\Omega$ of a bounded Sobolev
function, still denoted by 
$t\mapsto\psi(t)$, with $\psi(t)\in H^{1}(\Omega)\cap L^\infty(\Om)$. 

We assume also that $\psi\in
W^{1,\infty}(0,T;H^{1}(\Omega))\cap L^\infty({0,T};L^\infty(\Omega))$. Thus, the time
derivative $t\mapsto \dot{\psi}(t)$ belongs to the space
$L^\infty(0,T;H^{1}(\Omega))$ and its spatial gradient 
$t\mapsto D\dot{\psi}(t)$
belongs to the space $L^\infty(0,T;L^2(\Omega;\R ^2))$.
\end{subsection}

\begin{subsection}{ The external loads.}  We are interested in the case of
time-dependent dead loads, in which the density, 
$g\colon[0,T]\times\partial_N\Om\to\R $, of the applied surface force per unit area in the
reference configuration does not depend on the displacement $u$. We assume that 
the function $t\mapsto g(t,\cdot)$ belongs to
$W^{1,\infty}(0,T;L^2(\partial_N\Om,\hh))$, with time derivative denoted by $t\mapsto
\dot g(t,\cdot)$. The 
associated potential, for a displacement
$u$, is given by
$$
\GG(t)(u):=\int_{\partial_N\Om } g(t,x)u(x)d\hh \,.
$$
Moreover, assume that for every $t\in[0,T]$ the support of $g(t,\cdot)$ does not intersect the set $\Gamma$.

\end{subsection}
\begin{subsection}{The admissible displacements and their total energy}

For every $t\in [0,T]$, the set $AD(\psi(t),\sg)$ of admissible displacements in
$\Om$ with finite energy, corresponding to the crack $\G(\sg)$ and to the boundary
data
$\psi(t)$ 
 is given by
$$
AD(\psi(t),\sg):=\{u\in H^1(\Om\setminus\G(\sg)):u=\psi(t)\hbox{ on
}\partial_D\Om \}\,,
$$
where the last equality refers to the traces of $u$ and $\psi(t)$ on
$\partial_D\Om$. The total energy of a configuration $(u,\sg)$ with $u\in
AD(\psi(t),\sg)$ is given by
\begin{equation*}
\E(t)(u,\sg):=
\W(D u)+\sg %-\FF(t)(u)
-\GG(t)(u)\,.
\end{equation*}
Note that it does not depend on
the particular extension $\psi(t)$ chosen, but only on its value on the
Dirichlet part of the boundary.

\end{subsection}

\begin{subsection}{Moving to a fixed domain}
Let $H^1_{\partial_D\Om}(\Om\setminus\G(\sg))$ denote the space of functions $u\in
H^1(\Om\setminus\G(\sg))$ whose trace on $\partial_D\Om$ is zero. We may consider the
energy as a functional defined on 
$H^1_{\partial_D\Om}(\Om\setminus\G(\sg))$ by simply writing $\tilde u= u+\psi(t)$ with $\tilde u\in AD(\psi(t),\sg)$ and $u\in H^1_{\partial_D\Om}(\Om\setminus\G(\sg))$. 
Still the domain of the functional would depend on
$\sg$. To transform it into a functional defined on a fixed domain we consider the
following change of variables.

For $\sg\in[\sg_0,\overline\sg]$, let
$\Phi(\cdot,\sg)=\Phi_\sg(\cdot):\Om\to\Om$ be a
diffeomorphism which coincides with the identity near the boundary of $\Om$,
leaves invariant both $\Om^+$ and $\Om^-$ and transforms $\G(\sg)$ into
$\G(\sg_0)$. 
Let
$\Psi(\cdot,\sg)=\Psi_\sg(\cdot):=
\Phi^{-1}(\cdot,\sg)\colon\Om\to\Om$. 
Then
$$
\int_{\Om\setminus\G(\sg)}\!\!\!|D u+D\psi(t)|^2\!dx=
\int_{\Om\setminus\G(\sg_0)}\!\!\!|D u(\Psi_\sg(y))+
D\psi(t)(\Psi_\sg(y))|^2{\rm det}D\Psi_\sg(y)dy\,.
$$
For $u\in H^1_{\partial_D\Om}(\Om\setminus\G(\sg))$ define
$v(y,\sg):=u(\Psi_\sg(y))$ and let 
$\tilde \psi(t)(y,\sg):=\psi(t)(\Psi_\sg(y))$. 
With these notations
$$
\int_{\Om\setminus\G(\sg)}\!\!\!|D
u+D\psi(t)|^2\!dx\!=\!\!
\int_{\Om\setminus\G(\sg_0)}\!\!\!|((D\Psi_\sg)^T)^{-1}\!(y)(D v(y,\sg)+D
\tilde \psi(t)(y,\sg))|^2{\rm det}D\Psi_\sg(y)dy,
$$
and the last integral can be written also in the form
$$
\int_{\Om\setminus\G(\sg_0)}
\sum_{i,j\in\{1,2\}}a_{ij}(y,\sg)D_j(v(y,\sg)+\tilde\psi(t)(y,\sg))
D_i(v(y,\sg)+\tilde\psi(t)(y,\sg))dy\,,
$$
with the coefficients $a_{ij}$ given by the change of variables. 

Define $A(\sg):=(a_{ij}(\sg))_{ij}$ and note that $a_{ij}(\sg)\in C(\overline\Om)$, and 
$a_{ij}(\sg)=a_{ji}(\sg)$, for every $\sg\in[\sg_0,\overline\sg]$, and every $i,j$. 

We may assume that $0<c<\|{\rm det}D\Phi_\sg\|_\infty<C$ independently of $\sg\in[\sg_0,\overline\sg]$, where $\|\cdot\|_\infty$ denotes the $L^\infty$-norm on $\Om$.
Since $\Gamma$ is  of class $C^3$, we may also choose $\Phi(\cdot,\sg)$ (and hence $\Psi(\cdot,\sg)$) to depend regularly on $\sg$ in such a way that, as functions of $\sg$, the coefficients $a_{ij}$ be of class $C^2$ on $[\sg_0,\overline\sg]$, uniformly in $\overline\Om$. 
In particular, we shall use the fact that there exist positive constants $\lambda, \Lambda, \Lambda',\,L\,, L'>0$ independent of $\sg$, such that
\begin{equation}\label{ell}
(A(\sg)\xi|\xi)\geq\lambda|\xi|^2\qquad\forall\,\xi\in \R ^2\,,\quad\forall\,x\in\overline\Om\,,
\end{equation}
where 
$(\cdot|\cdot)$ denotes the scalar product in $\R ^2$,
\begin{eqnarray}
&&\|(A(\sg)\xi|\eta)\|_\infty\leq \Lambda |\xi|\,|\eta|\qquad\quad\forall\,\xi,\eta\in \R ^2\,,\label{lima}\\
&&\|(\partial_\sg A(\sg)\xi|\eta)\|_\infty\leq\Lambda' |\xi|\,|\eta|\qquad\forall\,\xi,\eta\in \R ^2\,,\label{limapr}\\
&&\|a_{ij}(\sg')-a_{ij}(\sg'')\|_\infty\leq
L|\sg'-\sg''|\qquad\hbox{and}\label{L}\\
&&\|\partial_\sg a_{ij}(\sg')-\partial_\sg a_{ij}(\sg'')\|_\infty\leq
L'|\sg'-\sg''|\label{L'}
\end{eqnarray}
for every $\sg',\sg''\in[\sg_0,\overline\sg]$ and $i,j=1,2$.

Note that, since
$\Psi_\sg$ coincides with the identity near the boundary of $\Om$, this change of
variables does not have any effect on $\GG$:  
$$
\GG(t)(u+\psi(t))=\GG(t)(v+\tilde\psi(t))\,.
$$
Moreover, we can neglect the dependence
of $\tilde \psi$ on $\sg$ since, for every 
$\sg\in [\sg_0,\overline\sg]$,
$\Psi_\sg$
coincides with the identity near the boundary of $\Om$, and we may assume that the
support of $\psi$ is included in the set where, for every 
$\sg\in [\sg_0,\overline\sg]$, $\Psi_\sg$ is the identity.
Hence the change of variables influences only the bilinear
term in $v$.

For brevity of notation, let 
$$
V:=H_{\partial_D\Om}^1(\Om\setminus \G(\sg_0))\,.
$$ 
On $V$
we consider the norm $\|\cdot\|_V$ defined by
$\|v\|_V:=\|Dv\|_2$, and the scalar product $(v,w)_V:=(Dv,Dw)$, where $\|\cdot\|_2$ and $(\cdot,\cdot)$ denote the norm and, respectively, the scalar product in $L^2(\Om)$ or
$L^2(\Om\setminus \G(\sg_0));\R ^2)$, depending on the context. 
Let
$V'$ denote its dual space and let $\langle \cdot,\cdot\rangle$ denote the duality
pairing between $V'$ and $V$. 

For every $t\in[0,T]$, $v\in V$, and $\sg\in [\sg_0,\overline\sg]$ 
define
\begin{eqnarray*}
&&\F(t,v,\sg):=\\
&&\ =\int_{\Om\setminus\G(\sg_0)}\!\!\!
\sum_{\ i,j\in\{1,2\}}\!\!\!
a_{ij}(\sg)D_j(v+\tilde\psi(t))D_i(v+\tilde\psi(t))dx+\sg%-\\
-\!\!\int_{\partial_N\Om} g(t)(v+\tilde\psi(t))d\hh.
\end{eqnarray*}
Then the functional
$\F$ can be also written as
\begin{eqnarray*}
\F(t,v,\sg):= \int_{\Om\setminus\G(\sg_0)} (A(\sg)D v| D v)\,dx + 
2\int_{\Om\setminus\G(\sg_0)}( D\psi(t)| Dv )\,dx-\\%-\int_\Om f(t,\sg)v\, dx-\\
 -\int_{\partial_N\Om} g(t)v\,d\hh 
+ \sg +
\int_\Om |D\psi(t)|^2dx-\int_{\partial_N\Om} g(t)\psi(t)\,d\hh,
\end{eqnarray*}
or
\begin{eqnarray*}
&&\!\!\! \!\!\! \F(t,v,\sg):=\\
&& =(A(\sg)D v, D v)+ 
2( D\psi(t), Dv ) -(g(t),v)_{\partial_N\Om} + \sg +
\|D\psi(t)\|_2^2-(g(t),\psi(t))_{\partial_N\Om}
\end{eqnarray*}
where
$(\cdot,\cdot)_{\partial_N\Om}$ denotes the scalar product in
$L^2(\partial_N\Om,\hh)$. 
Hence the elastic energy 
becomes 
$\F^{el}(t,v,\sg):= \F(t,v,\sg)-\sg$,  
and there exist four positive constants $\lambda_\F$, $\Lambda_\F$, $\mu_\F$, and $M_\F$, 
independent of $t$ and $\sg$, 
such that for every $t\in [0,T]$ and every $\sg\in[\sg_0,\overline\sg]$ 
\begin{eqnarray*}
&&\displaystyle \F^{el}(t,v,\sg)\geq \lambda_\F \| v\|_V^2 - \mu_\F \\
&&\displaystyle \F^{el}(t,v,\sg)\leq \Lambda_\F \|v\|_V^2  + M_\F,
\end{eqnarray*}
for every $v\in V$. 
Indeed, this follows from the uniform ellipticity of the bilinear part and 
standard estimates (on $\Om^+$ and $\Om^-$).
\end{subsection}

\begin{subsection}{Critical points of the energy}

For every $t\in [0,T]$ the function $\F(t,\cdot,\cdot)\colon V\times [\sg_0,\overline\sg]\to
\R $ is twice Fr\'echet partially differentiable with respect to $(v,\sg)$. The partial differential $\partial_v \F(t,v,\sg)$ belongs to $V'$, while 
the partial gradient $\grad_v \F(t,v,\sg)$ %, of $G$ with respect to $v$,
 is, by definition, the element of $V$ given by
$$
( \grad_v \F(t,v,\sg),w)_V = 
2(A(\sg)Dv,Dw) +2(\psi(t),w)_V-(g(t),w)_{\partial_N\Om}, 
$$ 
for every $w\in V$.
The partial differential 
$\partial_\sg \F(t,v,\sg)$ is given by
$$
\partial_\sg \F(t,v,\sg) = 
(\partial_\sg A(\sg)Dv,Dv) 
+ 1\,.
$$
For fixed $v\in V$ and $\sg\in[\sg_0,\overline\sg]$, we have that $\F(\cdot,v,\sg)\in W^{1,\infty}(0,T)$, with 
$$
\partial_t\F(t,v,\sg)=2(D\dot\psi(t),D v+D\psi(t)) -(\dot g(t),v+\psi(t))_{\partial_N\Om}-(g(t),\dot\psi(t))_{\partial_N\Om}\,.
$$
Note that by the regularity assumptions on $\psi$ and $g$ it follows also that
the map
$$
(t,v,\sg)\mapsto (\grad_v \F(t,v,\sg),\partial_\sg \F(t,v,\sg))
$$
is continuous from $(0,T)\times V\times (\sg_0,\overline\sg)$ into $V\times\R $.

The second order partial
differentials with respect to $(v,\sg)$ are given by
\begin{eqnarray*}
\langle\langle\partial^2_{(v,\sg)}\F(t,v,\sg)(w_1,\tau_1),(w_2,\tau_2)\rangle\rangle=(A(\sg)Dw_1,Dw_2)+(\partial_\sg
A(\sg)Dv,Dw_1)\tau_2+\\
+(\partial_\sg
A(\sg)Dv,Dw_2)\tau_1+
(\partial^2_{\sg\sg}A(\sg)Dv,Dv)\tau_1\tau_2\,,
\end{eqnarray*}
for every $(w_i,\tau_i)\in V\times\R $, $i=1,2$, where
$\langle\langle\cdot,\cdot\rangle\rangle$ denotes the duality product between
$V'\times\R $ and $V\times \R $.

Since, for fixed $t$ and $\sg$, the function $v\mapsto \F(t,v,\sg)$ is strictly convex, 
 it has a unique critical point $v_{t,\sg}$, and $v_{t,\sg}$ is a  minimum point. Also the function $u\mapsto \E(t)(u,\sg)$ is strictly convex and its critical
point is the unique 
 minimum point $u_{t,\sg}\in AD(\psi(t),\sg)$ of $u\mapsto
\E(t)(u,\sg)$. The function $u_{t,\sg}$  satisfies
$$
2\int_{\Om\setminus\G(\sg)}(D u_{t,\sg}|D w)dx=
\int_{\partial_N\Om}g(t,x)w\, d\hh\qquad\forall\,w\in
H^1_{\partial_D\Om}(\Om\setminus\G(\sg))\,.
$$

\begin{proposition}\llabel{p:minp}
For fixed $t\in [0,T]$ critical points of $\F(t,\cdot,\cdot)$ correspond to critical points of $\E(t)$ in the following sense: minimum points $v_{t,\sg}\in V$ of $v\mapsto \F(t,v,\sg)$ correspond by the change of variables to  minimum points $u_{t,\sg}\in AD(\psi(t),\sg)$ of $u\mapsto \E(t)(u,\sg)$. Moreover,  
$\partial_\sg \F(t,v_{t,\sg},\sg)=\partial_\sg  E(t,\sg)$, where $ E(t,\sg):=\E(t)(u_{t,\sg},\sg)$.
\end{proposition}

Before giving the proof we discuss some properties of the minimizers $u_{t,\sg}$.
The following result provides a useful characterization of the  ``singular" part
of the displacement $u_{t,\sg}$ near the tip $\gamma(\sg)$ of the crack.
For the
proof we refer to \cite{Gr-85}, \cite{Gr-92}.
\begin{proposition}\llabel{th:2.1}
Let $\sg\in[\sg_0,\overline\sg]$ and $u\in H^1(\Omega\setminus \G(\sg))$ be such that
\begin{equation}\llabel{L2}
\Delta u \in L^2(\Omega\setminus \G(\sg)) \quad \mbox{and} \quad 
\partial_\nu u=0 \; \mbox{ on }\G(\sg).
\end{equation}
Then there exists $\kk\in \R $ satisfying
\begin{equation}\llabel{kappa}
u-\kk\sqrt{r\frac{2}{\pi}} \sin\frac{\theta}{2}\in H^2(U\setminus \G(\sg)),
\end{equation} 
for every $U\subset\subset \Omega$ open. In (\ref{kappa}), 
$r(x):=|x-\gamma(\sg)|$ and $\theta(x)$ is the continuous function on 
$U\setminus \G(\sg)$ which coincides
with the counterclockwise oriented angle between $\dot{\gamma}(\sg)$ and 
$x-\gamma(\sg)$, and vanishes on the points of the form 
$x=\gamma(\sg)+h\dot{\gamma}(\sg)$ for $h>0$ sufficiently small.
\end{proposition}

The coefficient $\kk\sqrt{2/\pi}$ represents the {\em stress
intensity factor} associated to the displacement $u$ at the tip $\gamma(\sg)$.
We shall use its following characterization.

\begin{proposition}\llabel{p:kchar}
Let $\sg\in[\sg_0,\overline\sg]$, $u\in H^1(\Omega\setminus \G(\sg))$ satisfying
(\ref{L2}),  and let $\kk$ be defined by (\ref{kappa}). 
Then for every $\phi=(\phi_1,\phi_2)\in C^\infty_c(\Omega;\R ^2)$ we have
\begin{equation}\llabel{kchar}
\begin{array}{l}
{\displaystyle \kk^2 \phi(\gamma(\sg))\dot{\gamma}(\sg) = 
\int_{\Omega} \Big[ ((D_1 u)^2-(D_2 u)^2)(D_1 \phi_1-
D_2 \phi_2) +}\\
{\displaystyle +2D_1 u\,D_2 u(D_1\phi_2+D_2\phi_1)\Big]\, dx +
2\int_{\Omega} \Delta u (D_1 u\phi_1 + D_2 u\phi_2)\, dx.}
\end{array}
\end{equation}
\end{proposition}
\begin{proof}  For a complete proof see \cite[Proposition 2.2]{CD-DM-00}. The idea
is to consider 
$\eta>0$ such that
$\overline{B}(\gamma(\sg),\eta)\subset\Om$, to
integrate by parts:
\begin{equation}\llabel{kc0}
{\displaystyle\int_{\Omega\setminus B(\gamma(\sg),\eta)}
\Big[ ((D_1 u)^2-(D_2 u)^2)(D_1 \phi_1-D_2 \phi_2) 
+ 2D_1 u\,D_2 u(D_1\phi_2+D_2\phi_1)\Big]\, dx }
\end{equation}
and to pass to the limit as $\eta\to0$ using (\ref{L2}).
\end{proof}

\begin{proposition}\llabel{p:ensg}
The function $\sg\mapsto  E(t,\sg)$ is differentiable on
$[\sg_0,\overline\sg]$ and
\begin{equation}\llabel{ensg}
\partial_\sg E(t,\sg)= 1 - \kk_{t,\sg}^2\,,
\end{equation}
where $\kk_{t,\sg}\sqrt{\frac{2}{\pi}}$ is the  stress intensity factor associated to
$u_{t,\sg}$ at
$\gamma(\sg)$.
\end{proposition}
\begin{proof}
The proof follows the same arguments of the proof of \cite[Theorem~3.3]{CD-DM-00}.
To compute the partial derivative $\partial_\sg  E(t,\sg)$ we consider
a diffeomorphism similar to $\Phi_\sg$ and then use Proposition \ref{p:kchar}.
\end{proof}
Similar computations have been recently done in \cite{K-M} when the stored energy density $W$ is a polyconvex function with $W(A)=\infty$ for every matrix $A$ with $\det A\leq0$, and $\Gamma$ is a segment.

\begin{remark}\label{regmin} Fix $t_0\in{]0,T[}\,$. The map $\sg\mapsto v_{t_0,\sg}$ has the same regularity as $\sg\mapsto A(\sg)$, hence, under the regularity assumptions we made on $A(\sg)$, it is of class $C^2([\sg_0,\overline\sg])$. Since in this case we are not interested in the dependence on $t$, let us simplify the notation and set $v_\sg:=v_{t_0,\sg}$. Then  standard arguments for elliptic PDE's allow us to obtain that  for every $\sg*\in[\sg_0,\overline\sg]$ 
there exists $v'_{\sg*}\in V$ as strong limit in $V$ of the  difference quotient $\frac{v_{\sg}-v_{\sg*}}{\sg-\sg*}$, and the map $\sg\mapsto v'_\sg$ is continuous in the strong topology of $V$. The same arguments can be repeated to obtain that there exists $v''_{\sg^*}\in V$ as strong limit in $V$ of the difference quotient $\frac{v'_{\sg}-v'_{\sg*}}{\sg-\sg*}$ and that the map $\sg\mapsto v''_\sg$ is continuos with respect to the strong topology in $V$.
Note that $v'_\sg$ and $v''_\sg$ solve the following equations
\begin{eqnarray*}
&&(A(\sg)Dv'_\sg,Dw)+(\partial_\sg A(\sg)Dv_\sg,Dw)=0\qquad\qquad\qquad\qquad\qquad\quad\forall\, w\in V\,,\\
&&(A(\sg)Dv''_\sg,Dw)+2(\partial_\sg A(\sg)Dv'_\sg,Dw)+(\partial^2_\sg A(\sg)Dv_\sg,Dw)=0\quad\ \forall\, w\in V\,,
\end{eqnarray*}
respectively.
\end{remark}

\begin{proof}[Proof of Proposition~\ref{p:minp}] It follows from the change of
variables, Proposition~\ref{p:kchar}, and Proposition \ref{p:ensg}.
\end{proof}

\begin{remark}\label{r:posdef}
Fix $t_0\in{]0,T[}\,$. With the same notation as in Remark \ref{regmin} $v_\sg:=v_{t_0,\sg}$, note that the second order differential, $\partial^2_{(\sg,v)}\F(t_0,v_\sg,\sg)$, of $\F$ with respect to $(v,\sg)$ is strictly positive
definite if and only if the second order derivative of the function $\sg\mapsto
\F(t_0,v_\sg,\sg)$ is strictly positive, when both exist. Moreover,  by 
Proposition ~\ref{p:minp}, this is equivalent to the fact that the second order 
derivative  of $\sg\mapsto
 E(t_0,\sg)$ is strictly positive.
 
Indeed, %as for every $\sg\in[\sg_0,\overline\sg]$, $v_\sg$ is the minimum point of
as $\partial_v\F(t_0,v_\sg,\sg)=0$, and $\sg\mapsto v_\sg$ is, by Remark \ref{regmin}, a $C^2$-function, we have
$$
\langle\partial_\sg\partial_v
\F(t_0,v_\sg,\sg),w\rangle+\langle\partial^2_{vv}\F(t_0,v_\sg,\sg)
v'_\sg,w\rangle=0\quad\forall\,w\in V\,.
$$
Assume that
\begin{eqnarray*}
&&0<\frac{d}{d\sg}(\partial_\sg \F(t_0,v_\sg,\sg)+
\langle\partial_v \F(t_0,v_\sg,\sg), v'_\sg\rangle)=\\
&&\quad=\partial^2_{\sg\sg}\F(t_0,v_\sg,\sg)+
\langle\partial_\sg\partial_v \F(t_0,v_\sg,\sg), v'_\sg\rangle+
\langle\partial_v\partial_\sg \F(t_0,v_\sg,\sg),v'_\sg\rangle+\\
&&\quad+\langle\partial^2_{vv}\F(t_0,v_\sg,\sg)
v'_\sg,v'_\sg\rangle+\langle\partial_v
\F(t_0,v_\sg,\sg), v''_\sg\rangle\,.
\end{eqnarray*} 

Hence
$$
\partial^2_{\sg\sg} \F(t_0,v_\sg,\sg)+\langle\partial_v\partial_\sg
\F(t_0,v_\sg,\sg),v'_\sg\rangle>0\,,
$$
which implies that
$$
\partial^2_{\sg\sg}\F(t_0,v_\sg,\sg)>\langle\partial^2_{vv}\F(t_0,v_\sg,\sg)
v'_\sg,v'_\sg\rangle
$$
(recall that in our case $\langle\partial_\sg\partial_v
\F,w\rangle=\langle\partial_v\partial_\sg \F,w\rangle$).

Therefore 
\begin{eqnarray*}
&&\langle\langle\partial^2_{(v,\sg)}\F(t_0,v_\sg,\sg)(w,\tau),(w,\tau)\rangle\rangle=\\
&&\quad= \partial^2_{\sg\sg}\F(t_0,v_\sg,\sg)\tau^2+2
\langle\partial_\sg\partial_v\F(t_0,v_\sg,\sg),w\rangle\tau+\langle
\partial^2_{vv}\F(t_0,v_\sg,\sg)w,w\rangle>\\
&&\quad >
\langle\partial^2_{vv}\F(t_0,v_\sg,\sg)v'_\sg, v'_\sg\rangle\tau^2-
2\langle\partial^2_{vv}\F(t_0,v_\sg,\sg)v'_\sg,w\rangle\tau+\\
&&\qquad+\langle\partial^2_{vv}\F(t_0,v_\sg,\sg)w,w\rangle\geq 0\,,
\end{eqnarray*} 
which shows that
$\partial^2_{(\sg,v)}\F(t_0,v_\sg,\sg)$ is strictly positive definite.

It is also easy to see that if $\partial^2_{(\sg,v)}\F(t_0,v_\sg,\sg)$ is strictly
positive definite then the second order derivative of the function $\sg\mapsto \F(t_0,v_\sg,\sg)$ is strictly positive.
\end{remark}

\end{subsection}
\end{section}

\begin{section}{Irreversible quasistatic evolution}

Given an initial crack length $\sg_0>0$, and an initial value, $u_0$, of the
displacement, such that the initial configuration is in equilibrium, we want to study
a quasistatic evolution of  configurations $(u,\sg)$ which starts from $(u_0,\sg_0)$. 
We are
interested in the evolution until the crack length reaches the value 
$\sg_1$. We cannot avoid the solution to have jumps (even at $t=0$) to
configurations with crack lengths larger than $\sg_1$; if this is the case, then the
boundary data are
not compatible with a progressive crack growth on the interval
$[\sg_0,\sg_1]$.
\begin{definition}\llabel{d:qe}
The {\em irreversible quasistatic evolution problem} consists in finding a
left-continuous map
$t\mapsto(u(t),\sg(t))$, where $\sg(t)$ represents the length of the crack up to
time $t$, and the displacement $u(t)$  belongs to $AD(\psi(t),\sg(t))$, which
satisfies the following three conditions:
\begin{itemize}
\item[(a)] {\em local unilateral stability:} for every $t$
\begin{eqnarray}\llabel{ls1}
&& \E(t)(u(t),\sg(t))\leq \E(t)(u,\sg(t)) \qquad\forall u\in AD(\psi(t),\sg(t)) \\
\llabel{ls2}
&& \partial_\sg E(t,\sg(t))\geq 0,
\end{eqnarray} 
where $ E(t,\sg)$ is defined in Proposition \ref{p:minp};
\item[(b)] {\em irreversibility:} the map $t\mapsto \sg(t)$ is increasing;
\item[(c)] {\em energy inequality:} for every $0\leq s<t$ we have
\begin{eqnarray*}
&& \E(t)(u(t),\sg(t)) \leq \E(s)(u(s),\sg(s))+ \\
&&\qquad\qquad\qquad\qquad+\int_s^t\Big(2\int_{\Omega\setminus\G(\sg(\tau))}
(Du(\tau)|D\dot{\psi}(\tau)) dx-\int_{\partial_N\Om}\dot g(\tau)u(\tau)d\hh\Big)d\tau ,
\end{eqnarray*}
\end{itemize} 
\end{definition}

In terms of the functional $\F$, 
the  irreversible quasistatic evolution problem consists in finding a
left-continuous function
$t\mapsto (v(t),\sg(t))$ %from a time interval $[0,T_{max}]$ into $V\times\R $ 
which satisfies the
following three conditions:
\begin{itemize}
\item[($\rm a_\F$)] local unilateral stability: for every $t$ %\in[0,T_{max}]$,
\begin{equation}
\begin{cases}
\grad_v \F(t,v(t),\sg(t))=0, \\
\partial_\sg \F(t,v(t),\sg(t))\geq 0;
\end{cases}
\end{equation}
\item[($\rm  b_\F$)] irreversibility: the map $t\mapsto \sg(t)$ is increasing;
\item[($\rm c_\F$)] energy inequality: for every $0\leq s<t$ %<T_{max}$ 
we have
\begin{equation*}
\F(t,v(t),\sg(t))\leq \F(s,v(s),\sg(s)) + \int_s^t\partial_t
\F(\tau,v(\tau),\sg(\tau))d\tau\,.
\end{equation*}
\end{itemize} 
A solution, $t\mapsto (v(t),\sg(t))$, to this problem is called an {\em irreversible quasistatic evolution for $\F$}.% on the interval} $[0,T_{max}]$.

Let us remark that, by the very construction of the functional $\F$, an evolution for~$\F$ is well-defined only  for cracks whose length is less than or equal to $\overline\sg$.

In terms of an irreversible quasistatic evolution $t\mapsto(v(t),\sg(t))$ associated
to the functional $\F$, %on the time interval $[0,T_{max}]$,
 the Griffith's criterion can be
expressed as:
\begin{equation}\label{GrG}
\begin{cases}
\dot\sg(t)\geq0\\
\partial_\sg \F(t,v(t),\sg(t))\geq0\\
\partial_\sg \F(t,v(t),\sg(t))\dot\sg(t)=0
\end{cases}
\end{equation} 
for a.e.\ $t$. %\in[0,T_{max}]$.
Since the first two conditions are included in the definition of an
irreversible quasistatic evolution, it remains to prove the last one. 
\begin{proposition}\label{griffG}
Let $t\mapsto(v(t),\sg(t))$ be an %approximable 
irreversible quasistatic evolution for $\F$. % on $[0,T_{max}]$. 
Then for a.e.\ $t$ %\in[0,T_{max}]$
we have
$$
\partial_\sg \F(t,v(t),\sg(t))\dot\sg(t)=0\,.
$$
\end{proposition}
\begin{proof}
Since $\sg$ is increasing, $\dot\sg$ exists at a.e.\ $t$. %\in[0,T_{max}]$. 
Fix $t_0$ %\in [0,T_{max}]$
such that $\dot\sg(t_0)$ exists. As, given $\sg(t)$, the function $v(t)$ is determined
as the unique solution of $\grad_v \F(t,v,\sg(t))=0$, the hypotheses we made on $A(\sg)$ and on the
data $\psi$ and $g$ imply that $\dot v(t_0)$ exists, as strong
limit in $V$ of the difference quotient $\frac{v(t)-v(t_0)}{t-t_0}$.

As $\grad_v\F(t,v(t),\sg(t))=0$, from the energy inequality ${\rm (c_\F)}$ we deduce that for a.e.\
$t$ % \in[0,T_{max}]$ 
\begin{equation}\label{neg}
\partial_\sg \F(t,v(t),\sg(t))\dot\sg(t)\leq0\,.
\end{equation}
Since $\dot\sg(t)\geq0$ and $\partial_\sg \F(t,v(t),\sg(t))\geq0$, (\ref{neg}) implies the equality  to be proved. 

\end{proof}

Going back to the energy functional $\E$, the Griffith's criterion now reads
\begin{equation}\label{GrE}
\begin{cases}
\dot\sg(t)\geq0\\
1-\kappa^2(t)\geq0\\
(1-\kappa^2(t))\dot\sg(t)=0
\end{cases}
\end{equation} 
for a.e.\ $t$, where $\kappa(t)\sqrt{\frac{2}{\pi}}$ is the stress
intensity factor associated to the displacement $u(t)$ at the tip $\sg(t)$ (see
Proposition \ref{th:2.1}). Since by the change of variables we made,  
$\partial_\sg \F(t,v(t),\sg(t))=1-\kappa^2(t)$, the previous proposition shows that during an
irreversible quasistatic evolution the Griffith's criterion is satisfied. 
Note that this can be proved directly for $\E$, following, for instance, the lines of
\cite[Theorem 6.1]{DM-Toa-02}.

In the context of variational models for quasistatic crack propagation, the
 evolution of minimum energy configurations was studied (see, e.g. 
\cite{DM-T-02},
\cite{C},
\cite{F-L}, \cite{DM-F-T}) and existence results were proved in a very general
setting (see \cite{DM-F-T}). This kind of evolution is a solution to the following problem.
\begin{definition}\llabel{d:gsqe}
The {\em globally stable irreversible quasistatic evolution problem} consists in  
finding a solution to the irreversible quasistatic evolution problem 
which satisfies the {\em global stability condition:}
  for every $t$ %\in[0,T]$,
\begin{equation*}
\E(t)(u(t),\sg(t))\leq \E(t)(v,\sg)\qquad\forall\sg\geq\sg(t)\quad\forall v\in
AD(\psi(t),\sg)\,.
\end{equation*}
\end{definition}

During a globally stable irreversible quasistatic evolution the total energy is an absolutely continuous function of time and the energy inequality (c) becomes an equality.

However, a  solution to this
problem is not completely satisfactory since, in order to get the
global stability, we have to compare, at each time, the energy of a configuration with the energy of all admissible configurations with larger crack lengths. 
This is why we use another criterion of selection: among all irreversible quasistatic evolutions we choose the {\em approximable} ones, i.e. those that can be obtained as limits of solutions to a regularized evolution problem.

In this paper we consider the regularized problem given by 
a modified $\e$-gradient flow for the functional $\F$. 
Since we are interested in an irreversible crack growth for $\sg$ varying in the
interval $[\sg_0,\sg_1]$, we look for an increasing function $\sg(t)$. Hence, we consider the positive part of the derivative of
$\F$ with respect to
$\sg$.  Then, we  modify the evolution
law for the crack length  in such a way that it never reaches
$\overline\sg$. To this end we introduce a penalization factor $\lambda(\sg)$ that %will force the solution to remain in the interval $[\sg_0,\sg_2]$. This penalization factor
can be
any Lipschitz continuous function of $\sg$ which is equal to one for
$\sg\leq\sg_1$, is strictly positive for
$\sg_1<\sg<\overline\sg$, and is equal to
zero for
$\sg=\overline\sg$. For instance, let
\begin{equation}\llabel{lambda}
\lambda(\sg):=\frac{(\overline\sg-(\sg\lor\sg_1))^+}{\overline\sg-\sg_1}\,.
\end{equation}
In such a
way the evolution is the one given by the
$\e$-gradient flow, with the constraint that $\sg$ is increasing, on the
interval $[\sg_0,\sg_1]$ that we are interested in, and it is modified by
this artificial penalization term for $\sg>\sg_1$, so that we do not consider it
meaningful for $\sg>\sg_1$. 
\begin{definition}\label{weaksol}
A function $t\mapsto (v_\e(t),\sg_\e(t))$ is called a {\em solution to the
 initial value problem for the
 modified 
$\e$-gradient flow} %(\ref{defegradflow}) 
for the functional $\F$ on $[0,T]$
\begin{equation}\label{defegradflow}
\begin{cases}
\e\dot v_\e=-\grad_v \F(t,v_\e,\sg_\e)\\
\displaystyle\e\dot\sg_\e=(-\partial_\sg \F(t,v_\e,\sg_\e))^+\lambda(\sg_\e)\,,\\ 
v_\e(0)=u_0\\
\sg_\e(0)=\sg_0\,,
\end{cases}
\end{equation}
where $\lambda(\sg)$ is given by (\ref{lambda}),
if $v_\e\in C^1([0,T];V)$, $\sg_\e$ is a $C^1$ increasing function from $[0,T]$ into $[\sg_0,\overline\sg]$
and the first equation in (\ref{defegradflow}) is satisfied in the following sense
$$
(\e\dot v_\e,w)_V=-(\grad_v \F(t,v_\e,\sg_\e),w)_V\qquad\forall\,w\in V\quad\forall\, t\in[0,T]\,.
$$
\end{definition}
\noindent
Note that (\ref{defegradflow}) is a Cauchy problem for an ordinary differential equation in~$V\times\R $.

\begin{theorem}\label{thegradflow}
There exists a  solution $(v_\e,\sg_\e)$
to the initial value problem~(\ref{defegradflow})
with $\lambda(\sg)$ given by (\ref{lambda}), and the following energy estimate holds: for every $s$, $t\in[0,T]$ with $s<t$ 
\begin{eqnarray}
\e\int_s^t\|\dot
v_\e(\tau)\|^2_Vd\tau+\e\int_s^t\frac{|\dot\sg_\e(\tau)|^2}
{\lambda(\sg_\e(\tau))}\,d\tau+
\F(t,v_\e(t),\sg_\e(t))\leq\nonumber\\
\leq \F(s,v_\e(s),\sg_\e(s))+
\int_s^t\partial_t\F(\tau,v_\e(\tau),\sg_\e(\tau))\,d\tau.\label{epsenest}
\end{eqnarray}
 \end{theorem}
\begin{proof}
Taking into account the expressions of $\grad_v\F$ and $\partial_\sg \F$, the equations in 
(\ref{defegradflow}) can be written as
\begin{equation}\label{eqegradflow}
\begin{cases}
\e(\dot v_\e,w)_V=-2(A(\sg_\e)Dv_\e,Dw)-2(\psi(t),w)_V+(g(t),w)_{\partial_N\Om}\quad\forall\,w\in V\\
\e\dot\sg_\e=(-(\partial_\sg A(\sg_\e)Dv_\e,Dv_\e)-1)^+\lambda(\sg_\e)\,. 
\end{cases}
\end{equation}
Since the vector field defining the equation (\ref{eqegradflow}) depends on $t$ only through the boundary data $\psi$ and $g$, it is Lipschitz continuous in $t$. Moreover, for fixed $t$,
standard estimates show that it is Lipschitz continuous and bounded on the bounded subsets of $V\times\R $. Hence classical results on ODE's (see, e.g. \cite{D}) give the local existence and the uniqueness of the solution. Since there exist  $\alpha\in C([0,T])$ and $\beta>0$ such that
$$
(-\grad_v\F(t,v,\sg),v)_V+\sg(-\partial_\sg \F(t,v,\sg))^+\lambda(\sg)\leq \alpha(t)(\|v\|_V^2+\sg^2)+\beta
$$
for every $(v,\sg)\in V\times\R $,
the solution is defined on the whole interval $[0,T]$. 

The function $t\mapsto \F(t,v_\e(t),\sg_\e(t))$ is then Lipschitz continuous on $[0,T]$ with derivative  given for a.e.\ $t\in[0,T]$ by
\begin{eqnarray*}
\frac{d}{dt}\F(t,v_\e(t),\sg_\e(t))&=&\partial_t\F(t,v_\e(t),\sg_\e(t))+(\grad_v\F(t,v_\e(t),\sg_\e(t)),\dot v_\e(t))_V+\\
&&+\partial_\sg \F(t,v_\e(t),\sg_\e(t))\dot\sg_\e(t)\,.
\end{eqnarray*}
Taking into account the equations satisfied by $v_\e$ and $\sg_\e$, for every $s$, $t\in[0,T]$ with $s<t$ 
 we have
\begin{eqnarray*}
&& \F(t,v_\e(t),\sg_\e(t))-\F(s,v_\e(s),\sg_\e(s))=\\
&&\qquad\qquad=\int_s^t\left(\partial_t\F(\tau,v_\e(\tau),\sg_\e(\tau))-\e\|\dot v_\e(\tau)\|_V^2-\e\frac{(\dot\sg_\e(\tau))^2}{\lambda(\sg_\e(\tau))}\right)d\tau\,,
\end{eqnarray*}
which implies (\ref{epsenest}).
\end{proof}

\begin{remark}\label{sgsg1} Let $t\mapsto(v_\e(t),\sg_\e(t))$ be a solution to
(\ref{defegradflow}). %-(\ref{initegradflow}). 
Assume $\|v_\e(t)\|_V\leq M$ for some
positive constant
$M$ independent of
$t$ and $\e$. By (\ref{limapr}),
 $$
\e\dot\sg_\e(t)\leq (\Lambda'M^2+1)\lambda(\sg_\e(t))\leq C(\overline\sg-\sg_\e(t))^+\,,
$$
for some  constant $C>0$. 
By classical results on differential inequalities (see, e.g. \cite[Theorem I.6.1]{H})
it follows that for every $t\in [0,T]$
$$
\sg_\e(t)\leq\overline\sg-e^{-Ct/\e}(\overline\sg-\sg_0)\,,
$$
hence $\sg_\e$ never reaches $\overline\sg$. 
\end{remark}

Note that, since the evolution is constrained to cracks with lengths less than or equal to $\overline\sg$, Griffith's criterion is meaningful in this setting only until the length $\overline\sg$ is reached. As the penalization factor $\lambda(\sg)$ is strictly positive for $\sg<\overline\sg$, we may replace (\ref{GrG}) by
\begin{equation*} 
\begin{cases}
\dot\sg(t)\geq0\\
\partial_\sg \F(t,v(t),\sg(t))\lambda(\sg(t))\geq0\\
\partial_\sg\F(t,v(t),\sg(t))\dot\sg(t)=0\,.
\end{cases}
\end{equation*}
for a.e.\ $t\in [0,T]$. Therefore, also the second line in the local stability condition ${\rm (a_\F)}$ may be replaced by 
$\partial_\sg \F(t,v(t),\sg(t))\lambda(\sg(t))\geq0$.

We introduce now the following notion of evolution.
\begin{definition}\llabel{d:aqe} The {\em approximable irreversible quasistatic evolution problem} on the interval $[0,T]$
with initial data $(u_0,\sg_0)$ consists in finding a left-continuous map $t\mapsto(v(t),\sg(t))$ from $[0,T]$ into $V\times\R $ 
 %an irreversible  quasistatic evolution $t\mapsto (v(t),\sg(t))$  for $G$ %on some 
which satisfies the following conditions:
\begin{itemize}
\item[($\rm a'_\F$)] for every $t\in[0,T]$
\begin{eqnarray*}
&& \grad_v\F(t,v(t),\sg(t))=0\\
&& \partial_\sg \F(t,v(t),\sg(t))\lambda(\sg(t))\geq0\,;
\end{eqnarray*}
\item[($\rm b_\F$)] the map $t\mapsto\sg(t)$ is increasing;
\item[($\rm c_\F$)] for every $0\leq s<t\leq T$
$$
\F(t,v(t),\sg(t))\leq \F(s,v(s),\sg(s))+\int_s^t \partial_t \F(\tau,v(\tau),\sg(\tau))d\tau\,;
$$
\item[($\rm d_\F$)] the pair $(v(t),\sg(t))$ is the limit of a suitable solution 
$(v_\e(t),\sg_\e(t))$ 
of the {\em modified $\e$-gradient flow} for $\F$ 
with initial conditions $v_\e(0)=u_0$ and $\sg_\e(0)=\sg_0$, in the sense that for a.e. $t$ %\in[0,T_{max}]$
\begin{equation}\llabel{mm}
\begin{array}{l}
\sg_\e(t)\to\sg(t) ,\\
v_\e(t)\to v(t) \quad \mbox{strongly in }V.
\end{array}
\end{equation}
\end{itemize}
A solution $t\mapsto (v(t),\sg(t))$ to this problem is called an {\em approximable quasistatic evolution} for $\F$.
\end{definition}

We are now in a position to state the main result of this paper.
\begin{theorem}\label{main}
There exists a solution $t\mapsto(v(t),\sg(t))$ to the 
approximable irreversible quasistatic evolution problem with initial condition
$(u_0,\sg_0)$ on $[0,T]$. %on the maximal time interval on which $\sg(t)<\sg_2$.
\end{theorem}
\begin{remark} The fact that an approximable quasistatic
evolution starts from $(u_0,\sg_0)$ means only that for every
$\e>0$, $v_\e(0)=u_0$ and $\sg_\e(0)=\sg_0$.  We may always set
$(v(0),\sg(0)):=(u_0,\sg_0)$, but in general $v$ and $\sg$ are not continuous in $t=0$.
The only case in which $(u_0,\sg_0)$ is
the initial value for the evolution in a  ``classical" sense, is when $(u_0,\sg_0)$ is
the absolute minimum point of $\F(0,\cdot,\cdot)$. Indeed, in this case, by
semicontinuity and by the energy inequality ($\rm c_\F$), it is easy to see that
$t\mapsto \F(t,v(t),\sg(t))$ is continuous in $t=0$. 
\end{remark}

\begin{proof}[Proof of Theorem \ref{main}]
For $\e>0$ let $(v_\e,\sg_\e)$ be the solution of the modified $\e$-gradient flow with initial data $(u_0,\sg_0)$. Let $t\in[0,T]$.
The estimates we have on $\F$ together with (\ref{epsenest}) between $s=0$ and $t$ imply
$$
\lambda_\F\|v_\e(t)\|_V^2\leq \mu_\F+\F(0,u_0,\sg_0)+\int_0^t(a(\tau)\|v_\e(\tau)\|_V^2+b(\tau))d\tau
$$
for some functions $a,b\in L^\infty(0,T)$ which depend only on the data $\psi$ and $g$. Then, by Gronwall's Lemma, there exists a positive constant $C>0$ independent of $t$ and $\e$, whose value may change from line to line, such that
\begin{equation}\label{vH}
\|v_\e(t)\|_V\leq C \qquad\forall t\in[0,T]\,.
\end{equation}
By (\ref{epsenest}) we now get
\begin{eqnarray}
&\e\|\dot v_\e\|^2_{L^2(0,T;V)}\leq C\label{edotv2}\\
&\e\|\dot\sg_\e\|^2_{L^2(0,T)}\leq C\label{edotsg2}
\end{eqnarray}

Let $\e\to0$.
By Helly's Theorem, there
exists a subsequence, still denoted by $\e$, and an increasing function
$\sg\colon[0,T]\to[\sg_0,\overline\sg]$ such that
$$
\sg_\e(t)\to\sg(t)\qquad\hbox{for every }t\in[0,T]\,.
$$
The estimate (\ref{vH}) implies that there exists a function $v\in L^2(0,T;V)$ such that 
$$
v_\e\wto v \qquad \hbox{weakly in }L^2(0,T;V)\,,
$$
while, by (\ref{edotv2}), 
$$
\e\dot v_\e \to 0 \qquad \hbox{strongly in }L^2(0,T;V)\,.
$$
Hence
$$
\e(\dot v_\e(t),w)_V=(-\grad_v\F(t,v_\e(t),\sg_\e(t)),w)_V\to
(-\grad_v\F(t,v(t),\sg(t)),w)_V=0\,,
$$
for every $w\in V$ and for a.e. $t\in[0,T]$. It follows that 
$$
\int_0^T(A(\sg_\e(t))Dv_\e(t), Dv_\e(t))dt\to
\int_0^T(A(\sg(t))Dv(t),Dv(t))dt\,,
$$
which gives the strong convergence in $V$ of $v_\e(t)$ to $v(t)$  for a.e. $t\in[0,T]$.

By (\ref{edotsg2}),  
  $\e\dot\sg_\e(t)\to0$ for a.e. $t\in[0,T]$. Taking into account the equation satisfied by $\sg_\e$, we obtain that $(-\partial_\sg \F(t,v(t),\sg(t)))^+\lambda(\sg(t))=0$ for a.e.\ $t\in[0,T]$. 
  
  %Let $T_{max}:=\sup\{t\in[0,T]:\sg(t)<\sg_2\}$. Then 
When passing to the limit in (\ref{epsenest}), we neglect the terms containing the norms of the time derivatives of $v_\e$ and $\sg_\e$,  and thus get that for a.e.\ $s,\,t\in[0,T]$ with $s<t$ 
\begin{equation}\llabel{enestlim}
\F(t,v(t),\sg(t))\leq \F(s,v(s),\sg(s))+
\int_s^t\partial_t \F(\tau,v(\tau),\sg(\tau))d\tau\,.
\end{equation}
(By semicontinuity the estimate holds true for every $t\in[0,T]$.)

Since $\sg$ is increasing, for every $t\in[0,T]$ there exists the limit
$\sg^\ominus(t):=\lim_{s\to t-}\sg(s)$.  Let $v^\ominus(t)$ be the unique solution to
$\grad_v\F(t,v,\sg^\ominus)=0$. Then $v(s)\to v^\ominus(t)$ strongly in $V$ as
$s\to t-$, $\sg(t)=\sg^\ominus(t)$ and $v(t)=v^\ominus(t)$ for a.e.\ $t\in[0,T]$. 
By
construction, the map $t\mapsto (v^\ominus(t),\sg^\ominus(t))$ is left-continuous from $[0,T]$
into $V\times[\sg_0,\overline\sg]$. Moreover, 
$\partial_\sg \F(t,v^\ominus(t),\sg^\ominus(t))\lambda(\sg^\ominus(t))\geq0$ for every
$t\in[0,T]$. % with $\sg^\ominus(t)<\sg_2$. 
Let $s$, $t\in[0,T]$ with $s<t$, and let $s_n\to s-$, $t_n\to t-$ be such that
(\ref{enestlim}) holds for $s_n$ and $t_n$. Passing to the limit in (\ref{enestlim})
as
$n\to+\infty$ we obtain
\begin{equation*}
\F(t,v^\ominus(t),\sg^\ominus(t))\leq \F(s,v^\ominus(s),\sg^\ominus(s))+
\int_s^t\partial_t \F(\tau,v^\ominus(\tau),\sg^\ominus(\tau))d\tau\,,
\end{equation*}
so that we conclude that $(v^\ominus,\sg^\ominus)$ is an approximable quasistatic
evolution for $\F$ on $[0,T]$ which starts from $(u_0,\sg_0)$.
\end{proof}
\begin{remark}
From (\ref{enestlim}) we deduce that if $\overline t\in[0,T]$ is a discontinuity point of $t\mapsto \F(t,v(t),\sg(t))$  then 
$$
\lim_{t\to\overline t+} \F(t,v(t),\sg(t))\leq \F(\overline t,v(\overline t),\sg(\overline t))\,.
$$
Indeed, note that at every time $t$ the function $t\mapsto \sg(t)$ has a right limit. Let $\sg^\oplus(\overline t):=\lim_{t\to\overline t+}\sg(t)$, and let $v^\oplus(\overline t)$ be the solution to $\grad_v\F(\overline t,v, \sg^\oplus(\overline t))=0$. By the regularity assumptions made on the data, we have that $v(t)$ converges to $v^\oplus(\overline t)$ strongly in $V$, and hence, using (\ref{enestlim}), we obtain
$$
\lim_{t\to t+} \F(t,v(t),\sg(t))=\F(\overline t,v^\oplus(\overline t),\sg^\oplus(\overline t))\leq \F(\overline t,v(\overline t),\sg(\overline t))\,.
$$
\end{remark}
\end{section}

\begin{section}{Quasistatic evolution and the Implicit Function Theorem}

In this section we show that, under suitable regularity assumptions, the solution to the modified $\e$-gradient
flow converges to the continuous solution for the quasistatic evolution problem given by the
Implicit Function Theorem. 

\begin{theorem}\label{tfi}
Assume that in $(t^0,\sg^0)\in{[0,T[}\times[\sg_0,\sg_1[$ the following conditions are satisfied
\begin{eqnarray}
&&\partial_\sg E(t^0,\sg^0)=0\\
&&\partial^2_\sg E(t^0,\sg^0)>0\,.
\end{eqnarray}
Then there exists a  time interval $[t^0,t^1]$ and a unique Lipschitz continuous function
$\sg^0:[t^0,t^1]\to [\sg^0,\sg_1]$ such that
$$
\partial_\sg E(t,\sg^0(t))=0\qquad\forall\, t\in[t^0,t^1]\,.
$$
Moreover, if $(v_\e,\sg_\e)$ is the solution to the modified $\e$-gradient flow and the following two conditions are satisfied:  
\begin{eqnarray*}
&&\dot\sg_\e(t)>0\qquad\forall\,t\in[t^0,t^1]\\
&&\sg_\e(t^0)\to\sg^0\,,
\end{eqnarray*}
then $\sg_\e(t)\to\sg^0(t)$ and $E(t,\sg_\e(t))\to E(t,\sg^0(t))$ for every $t\in[t^0,t^1]$.
\end{theorem}

The first part of the theorem follows from the Implicit Function Theorem. As for the second part, let us remark that 
even if there are not at the moment general theorems guaranteeing the strict monotonicity of 
$\sg_\e$ during the approximation process, in many cases this will follow, for a suitable choice of the boundary data, from a symmetry argument.

We now prove the second part of the theorem in an equivalent form for the functional~$\F$. Indeed, since $\partial_\sg E(t,\sg)=\partial_\sg \F(t,v_{t,\sg},\sg)$ (see Proposition \ref{p:minp}), if the second order derivative $\partial^2_\sg E(t^0,\sg^0)>0$, then also $\frac{d}{d\sg}\partial_\sg \F(t^0,v_{t^0,\sg^0},\sg^0)>0$,  and this last condition is equivalent to the fact that the second order partial differential 
$\partial^2_{(v,\sg)}\F(t^0,v_{t^0,\sg^0},\sg^0)$ is strictly positive definite (see 
Remark~\ref{r:posdef}). 

\begin{theorem}\label{tfig}
Assume that in $(t^0,v^0,\sg^0)\in {[0,T[}\times V\times [\sg_0,\sg_1[$ the following
conditions are satisfied
\begin{equation*}
\begin{cases}
\grad_v\F(t^0,v^0,\sg^0)=0\,,\\
\partial_\sg \F(t^0,v^0,\sg^0)=0\,,
\end{cases}
\end{equation*}
and
the second order
differential,
$\partial^2_{(v,\sg)}\F(t^0,v^0,\sg^0)$, of
$\F$ with respect to $(v,\sg)$ is strictly positive definite, i.e. there exists
$\alpha>0$ such that
\begin{equation}\label{posdef0}
\langle\langle \partial^2_{(v,\sg)}\F(t^0,v^0,\sg^0)(w,\tau),(w,\tau)\rangle\rangle
\geq\alpha(\|w\|^2+|\tau|^2)\qquad\forall w\in V\quad\forall\tau\in\R \,.
\end{equation}

Let $(v_\e,\sg_\e)$ be the solution of the
modified
$\e$-gradient flow for
$\F$ given by Theorem~\ref{thegradflow} and assume that 
\begin{eqnarray}
&&v_\e(t^0)\to v^0\qquad\hbox{strongly in }V\quad\hbox{and}\nonumber\\
&&\sg_\e(t^0)\to\sg^0\qquad\hbox{as }\e\to0\,. \nonumber
\end{eqnarray}
Then there exist a time interval $[t^0,t^1]$ and a unique Lipschitz continuous function
$(v^0,\sg^0):[t^0,t^1]\to V\times  [\sg^0,\sg_1]$ such that
\begin{equation*}
\begin{cases}
\grad_v \F(t,v^0(t),\sg^0(t))=0\\
\partial_\sg \F(t,v^0(t),\sg^0(t))=0
\end{cases}
\end{equation*}
for every $t\in[t^0,t^1]$.
Assume that $\dot\sg_\e(t)>0$ and $\sg_\e(t)<\sg_1$ for every $t\in[t^0,t^1]$. Then 
$v_\e(t)\to v^0(t)$ strongly in $V$ and $\sg_\e(t)\to\sg^0(t)$ 
for every $t\in[t^0,t^1]$.
\end{theorem}

\begin{proof}
By our assumptions on the data, $\partial^2_{(v,\sg)}\F(t,v,\sg)$ (see
Subsection 2.7) is continuous with respect to $(t,v,\sg)\in [0,T]\times
V\times[\sg_0,\overline\sg]$. Moreover,  the function
$t\mapsto\partial_t\grad_v\F(t,v,\sg)$ belongs to $L^{\infty}(0,T;V)$, while
$\partial_t\partial_\sg \F(t,v,\sg)=0$. By the Implicit Function Theorem (see, e.g., \cite{K})
 applied in $(t^0,v^0,\sg^0)$ to 
\begin{eqnarray*}
\begin{cases}\grad_v \F(t,v,\sg)=0\\
\partial_\sg \F(t,v,\sg)=0\,,
\end{cases}
\end{eqnarray*}
it follows that there exist a time interval $[t^0,t^1]$ and a unique Lipschitz continuous
function
$(v^0,\sg^0)\colon [t^0,t^1]\to V{\times}[\sg^0,\sg_1)$ such that
\begin{equation}\llabel{00}
\begin{cases}
\grad_v \F(t,v^0(t),\sg^0(t))=0\\
\partial_\sg \F(t,v^0(t),\sg^0(t))=0
\end{cases}
\end{equation}
for every $t\in[t^0,t^1]$. By a compactness argument, changing eventually
the value of $\alpha$, we may assume that 
there exist $\alpha>0$ and $r>0$ such that for every 
$t\in[t^0,t^1]$,  
for every $v\in B_r(v^0(t))\subset V$, and for every $\sg\in (\sg^0(t)-r,\sg^0(t)+r)$
\begin{equation}\label{posdef1}
\langle\langle \partial^2_{(v,\sg)}\F(t,v,\sg)(w,\tau),(w,\tau)\rangle\rangle
\geq\alpha(\|w\|^2+|\tau|^2)\qquad\forall w\in V\quad\forall\tau\in\R \,.
\end{equation}
Restricting eventually the time interval, we have
$\sg^0(t)+r<\sg_1$ for every $t\in[t^0,t^1]$.

Let $0<r'<r$ be a number that we shall choose later. For every $\e>0$ small enough we have $\|v_\e(t^0)-v^0\|_V<r'$ and $|\sg_\e(t^0)-\sg^0|<r'$.
By continuity, there exists a time interval, depending on $\e$, on which these inequalities hold. Let $\tau_\e$ be the largest time 
such that for
$t<\tau_\e$,  
$\|v_\e(t)-v^0(t)\|_V<r'$ and $|\sg_\e(t)-\sg^0(t)|<r'$. Then $\sg_\e(t)<\sg_1$, hence
$\lambda(\sg_\e(t))=1$ for $t<\tau_\e$. 

We want to prove that $\tau_\e=t^1$. Assume by contradiction that $\tau_\e<t^1$.
Taking $v_\e(t)-v^0(t)$ as test function in the equation satisfied by $v_\e$,
multiplying by $\sg_\e(t)-\sg^0(t)$ the equation satisfied by $\sg_\e$, and taking
also into account (\ref{00}), we obtain
\begin{eqnarray*}
&&\frac{\e}{2}\frac{d}{dt}\|v_\e(t)-v^0(t)\|_V^2+
\frac{\e}{2}\frac{d}{dt}|\sg_\e(t)-\sg^0(t)|^2=\\
&&=-(\grad_v \F(t,v_\e(t),\sg_\e(t))-
\grad_v \F(t,v^0(t),\sg^0(t)),v_\e(t)-v^0(t))_V
+\\
&&\qquad+\big(-\partial_\sg \F(t,v_\e(t),\sg_\e(t))
+\partial_\sg \F(t,v^0(t),\sg^0(t))\big)(\sg_\e(t)-\sg^0(t))-\\
&&\qquad-\e(\dot
v^0(t),v_\e(t)-v^0(t))_V-\e\dot\sg^0(t)(\sg_\e(t)-\sg^0(t))\,.
\end{eqnarray*}
 Setting
$$
\zeta_\e(t):=\|v_\e(t)-v^0(t)\|_V^2+|\sg_\e(t)-\sg^0(t)|^2\,,
$$
from (\ref{posdef1}) it follows that 
\begin{eqnarray}
&&\!\!\!\!\!\frac{\e}{2}\dot\zeta_\e(t)\leq -\alpha\zeta_\e(t)
-\e(\dot v^0(t),v_\e(t)-v^0(t))_V-\e\dot
\sg^0(t)(\sg_\e(t)-\sg^0(t))\leq\nonumber\\
&&\leq-\alpha\zeta_\e(t)+\frac{\e}{2}\|\dot v^0(t)\|_V^2+\frac{\e}{2}\|v_\e(t)-v^0(t)\|_V^2+
\frac{\e}{2}|\dot\sg^0(t)|^2+
\frac{\e}{2}|\sg_\e(t)-\sg^0(t)|^2\leq\nonumber\\
&&\leq(-\alpha+\frac{\e}{2})\zeta_\e(t)+
\frac{\e}{2}\beta\qquad\forall t\in[t^0,\tau_\e)\,,\llabel{stimazeta}
\end{eqnarray}
where $\beta$ is an upper bound for 
$\|\dot v^0(t)\|_V^2+|\dot\sg^0(t)|^2$ on $[t^0,t^1]$.

Hence 
\begin{equation}\label{beta}
\zeta_\e(t)\leq
\Big(\zeta_\e(t^0)-\frac{\beta\e}{2\alpha-\e}\Big)
e^{(-\frac{2\alpha}{\e}+1)(t-t^0)}
+\frac{\beta\e}{2\alpha-\e}\qquad\forall\, t\in[t^0,\tau_\e)\,.
\end{equation}
Therefore, choosing now $r'$ small enough, from (\ref{beta}) we get that also 
$\|v_\e(\tau_\e)-v^0(\tau_\e)\|_V<r$ and
$|\sg_\e(\tau_\e)-\sg^0(\tau_\e)|<r$. By continuity, these inequalities hold also for some $t>\tau_\e$, which contradicts the maximality of $\tau_\e$, and so we deduce that $\tau_\e=t^1$. Hence (\ref{beta}) holds for
every
$t\in[t^0,t^1]$. Passing to the limit in (\ref{beta}) we get the conclusion.
\end{proof}

By the change of variables that defines the functional $\F$, and by the uniqueness of the regular evolution given by the Implicit Function Theorem, it follows that the regular evolution in Theorem~\ref{tfig} corresponds to the one in Theorem~\ref{tfi}. 

\begin{proof}[Proof of Theorem \ref{tfi} continued]
Let $(v_\e,\sg_\e)$ be the solution to the modified $\e$-gradient flow. By Theorem~\ref{tfig}, $\sg_\e(t)\to\sg^0(t)$ and $v_\e(t)\to v^0(t)$ strongly in $V$ for every $t\in[t^0,t^1]$. Since the function $v\mapsto\grad_v\F(t,v,\sg)$ is continuous from $V$ to $V$ with respect to the strong topology, it follows that
$$
\grad_v\F(t,v_\e(t),\sg_\e(t))\to\grad_v\F(t,v^0(t),\sg^0(t))=0\,.
$$ 

 Let $\overline v_\e(t)$ be the element of $V$ associated to $u_{t,\sg_\e(t)}$ by the change of variables.  As $\grad_v\F(t,\overline v_\e(t),\sg_\e(t))=0$ we deduce that $v_\e(t)-\overline v_\e(t)\to0$ strongly in $V$. This implies that 
$$
\F(t,\overline v_\e(t),\sg_\e(t))-\F(t,v_\e(t),\sg_\e(t))\to0\,, 
$$
On the other hand, 
\begin{eqnarray*}
&&\F(t,v_\e(t),\sg_\e(t))\to\F(t,v^0(t),\sg^0(t))=E(t,\sg^0(t))\\
&&\F(t,\overline v_\e(t),\sg_\e(t))= E(t,\sg_\e(t))\,,
\end{eqnarray*}
 so that we conclude that $E(t,\sg_\e(t))\to E(t,\sg^0(t))$ for every $t\in[t^0,t^1]$. 
\end{proof} 
\end{section}

\begin{section}{Monotonically increasing loadings}

In  this section we consider the setting proposed by 
Francfort and Marigo in \cite{FrMa-98} and compare the  evolution defined therein 
with a solution to the irreversible quasistatic evolution problem. 

Assume $\psi(t):=t\psi_0$, with
$\psi_0\in H^1(\Omega)$, and $g(t)=0$, and define
\begin{equation*}
E(\sg):=\min \{\|D u\|_2^2: u\in AD(\psi_0,\sg)\}.
\end{equation*}
Since $H^1(\Om\setminus\Gamma(\sg'))\subset H^1(\Om\setminus\Gamma(\sg''))$ for
$\sg'<\sg''$, we have that $E(\sg')\geq E(\sg'')$, so that the function
$\sg\mapsto E(\sg)$ is decreasing.

Following Definition~4.13 in \cite{FrMa-98}, we define a crack trajectory 
$t\mapsto \sg_{\scriptscriptstyle FM}(t)$ by the following three properties:
\begin{itemize}
\item[(i)] $t\mapsto \sg_{\scriptscriptstyle FM}(t)$ is increasing;
\item[(ii)] $t^2 E(\sg_{\scriptscriptstyle FM}(t))+\sg_{\scriptscriptstyle FM}(t)\leq
t^2 E(\sg) +\sg$ for every
$\sg\geq \sg_{\scriptscriptstyle FM}^-(t)$;
\item[(iii)] $t^2 E(\sg_{\scriptscriptstyle FM}(t))+\sg_{\scriptscriptstyle FM}(t)\leq
t^2E(\sg_{\scriptscriptstyle FM}(s)) + \sg_{\scriptscriptstyle FM}(s)$ for every
$s\leq t$.
\end{itemize}
The following result shows that if $\sg\to E(\sg)$ is concave in some
subinterval of 
$(\sg_0,\overline\sg)$ then $t\mapsto \sg_{\scriptscriptstyle FM}(t)$ is discontinuous.
\begin{proposition}\llabel{p:disc}
Let $t\mapsto \sg_{\scriptscriptstyle FM}(t)$ be a crack trajectory which satisfies
properties (i)--(iii) above. If there exists a subinterval
$(a,b)\subset (\sg_0,\overline\sg)$ where 
$\sg\mapsto E(\sg)$ is concave, then $\sg_{\scriptscriptstyle FM}(t)$ has some
discontinuity points.
\end{proposition}
\begin{proof} 
Let $t_0\geq 0$ be such that $\sg_{\scriptscriptstyle FM}(t_0)<a$. We first prove
that there exists
$t>t_0$  such that $\sg_{\scriptscriptstyle FM}(t)>a$.
Indeed, assume by contradiction that $\sg_{\scriptscriptstyle FM}(t)<a$ for every
$t>t_0$. Then conditions  (i), (ii) and the fact that $\sg\mapsto E(\sg)$ is
decreasing imply the following inequalities:
\begin{equation*}
t^2 E(b) + b \geq t^2 E(\sg_{\scriptscriptstyle FM}(t)) +
\sg_{\scriptscriptstyle FM}(t) \geq t^2 E(a) +
\sg_{\scriptscriptstyle FM}(t)
\geq t^2 E(a) +  \sg_{\scriptscriptstyle FM}(t_0).
\end{equation*}
In particular, we deduce that $t^2\leq b(E(a)-E(b))^{-1}$, which, up to
considering  $T$ large enough, represents a contradiction.

If $\sg_{\scriptscriptstyle FM}(t)\neq a$ 
for every $t\in[0,T]$ then $\sg_{\scriptscriptstyle FM}$ is discontinuous
and the proof is concluded. Otherwise, let
$\bar{t}$ be the first time such that $\sg_{\scriptscriptstyle FM}(\bar{t})=a$.
We claim that
\begin{equation}\llabel{disc1}
\sg_{\scriptscriptstyle FM}(t)=a \quad \mbox{for every } \bar{t}\leq t \leq t^*,
\end{equation}
where 
\begin{equation}\llabel{t*}
t^*:= \sqrt{\frac{b-a}{E(a)-E(b)}}.
\end{equation}
Indeed, fix $t\in(\bar{t},t^*)$ and assume by contradiction that 
$\sg_{\scriptscriptstyle FM}(t)\in
{]a,b]}$.  Then there exists $\alpha\in ]0,1]$ such that 
$\sg_{\scriptscriptstyle FM}(t)=\alpha
a +(1-\alpha)b$. By condition (ii) and the concavity of $\sg\mapsto E(\sg)$ we
have
\begin{equation*}
t^2 E(b) + b \geq t^2E(\sg_{\scriptscriptstyle FM}(t)) +
\sg_{\scriptscriptstyle FM}(t) \geq t^2\alpha E(a) + t^2(1-\alpha) E(b) +
\alpha
a+(1-\alpha)b,
\end{equation*}
that is 
\begin{equation}\llabel{disc2}
t^2 E(a) +a \leq t^2 E(b) + b.
\end{equation}
Therefore
\begin{equation}\llabel{disc3}
t^2 E(\sg_{\scriptscriptstyle FM}(t)) + \sg_{\scriptscriptstyle FM}(t) \geq
\alpha(t^2 E(a) + a) + (1-\alpha)(t^2 E(b) +b)
\geq t^2
 E(a) + a.
\end{equation}
Since (\ref{disc3}) is in 
contradiction with condition (iii), we deduce that $\sg_{\scriptscriptstyle
FM}(t)=a$.

Consider now the case $t=t^*$. Formula (\ref{disc2}) becomes the identity
\begin{equation*}
(t^*)^2 E(a) +a = (t^*)^2 E(b) + b.
\end{equation*}
Assume there exists $\alpha\in ]0,1]$ such that $\sg_{\scriptscriptstyle
FM}(t^*)=\alpha a + (1-\alpha)b$; then,  arguing as before, we obtain that
\begin{equation*}
(t^*)^2 E(\sg_{\scriptscriptstyle FM}(t^*)) + \sg_{\scriptscriptstyle FM}(t^*)\geq
(t^*)^2 E(a) +a,
\end{equation*}
which, by conditions (ii) and (iii), implies that $\sg_{\scriptscriptstyle
FM}(t^*)=a$.

To conclude, we prove that 
\begin{equation}\llabel{disc4}
\sg_{\scriptscriptstyle FM}(t)\geq b \quad \mbox{for }t> t^*.
\end{equation}
Indeed, let us fix $t>t^*$ and assume by contradiction that $\sg_{\scriptscriptstyle
FM}(t)<b$.  Then there exists $\alpha\in ]0,1]$ such that
$\sg_{\scriptscriptstyle FM}(t)=\alpha a + (1-\alpha)b$, and this fact
together with condition (ii) implies that
\begin{equation}
t^2\leq \frac{b-\sg_{\scriptscriptstyle FM}(t)}{E(\sg_{\scriptscriptstyle
FM}(t))-E(b)} \leq
\frac{\alpha(b-a)}{\alpha E(a) + (1-\alpha)E(b)- E(b)} = (t^*)^2,
\end{equation}
a contradiction. This fact concludes the proof, since we have shown that for 
$t\leq t^*$ $\sg_{\scriptscriptstyle FM}(t)=a$, while $\sg_{\scriptscriptstyle
FM}(t)\geq b$ for $t>t^*$.
\end{proof}

Let $(u(\cdot),\sg(\cdot))$ be an irreversible quasistatic evolution. 
Recalling that $u(t)$ is the minimum point of $\|D u\|_2$ on $AD(t\psi_0,\sg(t))$, we
have that $\|D u(t)\|_2^2=t^2 E(\sg(t))$.
We may now express conditions $(a)$, $(b)$ and $(c)$ of
Definition~\ref{d:qe} of an  irreversible quasistatic evolution, in terms of $\sg(t)$,
and, in the case of this particular choice of the data, we obtain:  
\begin{itemize}
\item[$(a')$] $1+t^2E'(\sg(t))\geq 0$ for every $t\geq 0$;
\item[$(b')$] the map $t\mapsto \sg(t)$ is
increasing;
\item[$(c')$] ${\displaystyle t^2 E(\sg(t))+\sg(t)\leq s^2 E(\sg(s)) + \sg(s) +
2\int_s^t\tau  E(\sg(\tau))\, d\tau}$, for every ${0\leq s<t}$, 
\end{itemize}
where $E'(\sg(t))$ denotes the derivative of $E$ with respect to $\sg$ computed at 
$\sg(t)$. 

Since $E(\sg(\tau))\leq E(\sg(s))$ for $\tau\in[s,t]$,
condition $(c')$ implies  condition (iii).
\begin{remark}\llabel{r:Esg4}
Let $t\mapsto \sg(t)$ be a left-continuous map on $[0,T]$ which satisfies 
condition $(c')$ and define
\begin{equation*}
\dot{\sg}^{\ominus}(t):= \limsup_{s\to t^-} \frac{\sg(t)-\sg(s)}{t-s}. 
\end{equation*}
Then 
\begin{equation}\llabel{Esg4}
(1+t^2E'(\sg(t)))\dot{\sg}^{\ominus}(t)\leq 0
\end{equation}
for every $t\in[0,T]$. Indeed, let $t_k\nearrow t$ be such that
\begin{equation*}
\lim_{k\to\infty} \frac{\sg(t)-\sg(t_k)}{t-t_k}=\dot{\sg}^{\ominus}(t).
\end{equation*}
Then condition $(c')$ between $t_k$ and $t$ can be written as
\begin{equation*}
(t^2-t_k^2)E(\sg(t)) + t_k^2(E(\sg(t))-E(\sg(t_k))) + \sg(t)-\sg(t_k)\leq 
2\int_{t_k}^t\tau E(\sg(\tau))\, d\tau,
\end{equation*}
and (\ref{Esg4}) follows dividing by $t-t_k$ and letting $k\to+\infty$. 
\end{remark}

\begin{remark}\llabel{r:disc}
Let $t\mapsto \sg(t)$ be a left-continuous map on $[0,T]$ which satisfies 
conditions $(a')$, $(b')$, and $(c')$, and 
let $t\geq 0$ be such that $\dot{\sg}^{\ominus}(t)>0$. Then, by Remark \ref{r:Esg4}
and conditions $(a')$ and $(b')$, it follows that
\begin{equation*}
E'(\sg(t))=\frac{d}{d\sg}E(\sg)|_{\sg=\sg(t)}= - \frac{1}{t^2}\,,
\end{equation*}
which implies that $\sg(t)$ does not belong to the concavity intervals of $E(\sg)$, 
since $t\mapsto \sg(t)$ is increasing, and $t\mapsto E'(\sg(t))$ would be decreasing,
while the right-hand side is increasing. More precisely, if there exists an interval
$(a,b)\subset[\sg_0,\overline\sg]$ such that $\sg\mapsto E'(\sg)$ is strictly  decreasing 
on $(a,b)$ and there exists
$t_0\geq 0$ such that
$\dot{\sg}(t_0)>0$ (or $\dot{\sg}^\ominus(t_0)>0$) and $\sg(t_0)\in (a,b)$ then we
reach a contradiction. Indeed, let $t>t_0$ be such that $\sg(t)\in
(a,b)$. By $(b')$, $\sg(t)>\sg(t_0)$, and by $(a')$ and our assumption on $E'(\sg)$,
we  get
\begin{equation*}
-\frac{1}{t^2}\leq E'(\sg(t))<E'(\sg(t_0))=-\frac{1}{t_0^2} < -\frac{1}{t^2},
\end{equation*}
a contradiction.
\end{remark}

In order to specify better the monotonicity needed in the above remarks we introduce
the following notion. We say that $t_0$ is a {\em local left-constancy point} for
$\sg$ if there exists
$\e>0$ such that $\sg$ is constant on the interval $[t_0-\e,t_0]$.

\begin{proposition}\llabel{p:Econc}
Let $\sg\colon[0,T]\to[\sg_0,\overline\sg[$ be a left-continuous map which satisfies conditions
$(a')$, $(b')$, and $(c')$, and let $t_0\geq 0$. If
\begin{itemize}
\item[(1)] $t_0$ is not a local left-constancy point for $\sg$ and 
\item[(2)] there exists $(a,b)\subset[\sg_0,\overline\sg[$ such that $E'(\sg)$ is strictly 
decreasing on $(a,b)$
\end{itemize}
then $\sg(t_0)\notin (a,b)$.
\end{proposition}
\begin{proof}
If $t_0$ is not a local left-constancy point for $\sg$, then, given $\e>0$, there are
$t^1_\e,t^2_\e\in [t_0-\e,t_0]$ such that 
$\sg(t^1_\e)\neq \sg(t^2_\e)$. Therefore, there exists $t_\e\in {[t_0-\e,t_0[}$ such
that $\dot{\sg}^\ominus(t_\e)>0$. Then (\ref{Esg4}) together with $(a')$ imply that
$1+t^2_\e E'(\sg(t_\e))=0$.
By Remark \ref{r:disc}, $\sg(t_\e)\notin (a,b)$ and we conclude by passing to the
limit as
$\e\to0$ (since $\sg$ is left-continuous).
\end{proof}
\begin{proposition}\llabel{p:Econv}
Let $\sg:[0,T]\to[\sg_0,\overline\sg[$ be a left-continuous map which satisfies conditions
$(a')$, $(b')$, and $(c')$.
Assume that $E(\sg)$ is convex on $(a,b)\subset[\sg_0,\overline\sg]$. Then $\sg(t)$ is 
continuous at every $t$ with $\sg(t)\in (a,b)$.
\end{proposition}
\begin{proof}
Assume by contradiction that $\sg(t)<\sg(t^+)$. Then condition $(c')$ and condition 
$(a')$ imply 
\begin{equation*}
\frac{E(\sg(t^+)-E(\sg(t)))}{\sg(t^+)-\sg(t)}\leq -\frac{1}{t^2}\leq E'(\sg(t)),
\end{equation*}
a contradiction.
\end{proof}
\end{section}
\begin{section}{Concavity and convexity intervals for the energy
functional}\label{enconcave}

We consider the energy functional
$$
\sg\mapsto E(\sg):=\min\{\|D u\|_2^2:u\in AD(\psi,\sg)\},
$$
and construct an explicit example of $\Om$ and $\psi$ for which $E(\sg)$ is concave on some subinterval. 
 Let $B_{-2}$ denote the ball of radius $1$ centred in $(-2,0)$, let $B_2$ denote the
ball of radius $1$ centred in $(2,0)$, and let $\G:= [-3,3]\times \{0\}$.
\begin{figure}[ht]
\begin{center}
\psfrag{0}{$0$}
\psfrag{-2}{$-2$}
\psfrag{2}{$2$}
\psfrag{G}{$\Gamma$}
\psfrag{B-}{$B_{-2}$}
\psfrag{B}{$B_2$}
\psfrag{T}{$T_{\varepsilon}$}
\psfrag{theta}{$\theta$}
\psfrag{ttheta}{$\tilde{\theta}$}
\includegraphics[width=300pt]{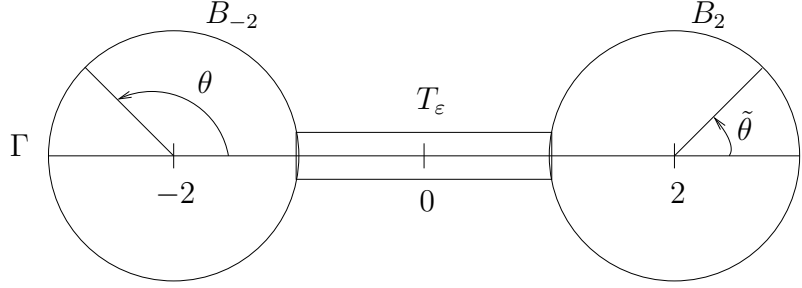}
\caption{\em The set $\Omega_\varepsilon$.}
\end{center}
\end{figure}

For $\e>0$ let
\begin{equation*}
T_\e := {]-2+\cos\e,2-\cos \e[}\times {]-\sin \e, \sin \e[}\, , \qquad
\Omega_\e:= B_{-2}\cup T_\e \cup B_{2}.
\end{equation*}
Further, for every $\sg\in[-3,3]$ let
\begin{equation*}
\G(\sg):=[-3,\sg]\times \{0\}.
\end{equation*}
Let $(\rho,\theta)$ and $(\tilde\rho,\tilde \theta)$ be polar coordinates around
$(-2,0)$ and $(2,0)$, respectively, where the functions  $\theta$ and $\tilde\theta$
are chosen, as in Proposition \ref{th:2.1}, such that $\theta(x_1,x_2)\to-\pi$ if
$x_2\to0-$ and $x_1<-2$, $\theta(x_1,x_2)\to\pi$ if $x_2\to0+$ and $x_1<-2$, and,
analogously, $\tilde\theta(x_1,x_2)\to-\pi$ if $x_2\to0-$ and $x_1<2$,
$\tilde\theta(x_1,x_2)\to\pi$ if $x_2\to0+$ and $x_1<2$.

On $\partial\Om_\e$ we define the boundary data $\psi_\e$ as follows:
\begin{equation}\llabel{psi}
\psi_\e(x):=
\begin{cases}
\sin\frac{\theta(x)}{2} &\mbox{ on } 
\big(\partial B_{-2}\cap\partial\Omega_\e\big)\setminus \G(\sg),\\
\sin\frac{\tilde{\theta}(x)}{2} &\mbox{ on } 
\big(\partial B_{2}\cap\partial\Omega_\e\big)\setminus \G(\sg),\\
\sin\frac{\e}{2} &\mbox{ on } {]-2+\cos\e,0[}\times\{\sin \e\},\\
-\sin\frac{\e}{2} &\mbox{ on } {]-2+\cos\e,0[}\times\{-\sin \e\},\\
\sin\frac{\e}{2} + \frac{x_1}{2-\cos\e}\big(\cos\frac{\e}{2} -
\sin\frac{\e}{2}\big)&\mbox{ on }  {[0,2-\cos \e[}\times\{\sin \e\},\\
-\sin\frac{\e}{2} + \frac{x_1}{2-\cos\e}\big(\sin\frac{\e}{2} - \cos\frac{\e}{2}
\big)&\mbox{ on }  {[0,2-\cos \e[}\times\{-\sin \e\}.
\end{cases}
\end{equation}
\begin{figure}[ht]
\begin{center}
\psfrag{e}{\scriptsize${\varepsilon}$}
\psfrag{sint2}{\footnotesize $\sin\frac{\theta}{2}$}
\psfrag{sine2}{\footnotesize$\sin\frac{\varepsilon}{2}$}
\psfrag{-sine2}{\footnotesize$-\sin\frac{\varepsilon}{2}$}
\psfrag{cose2}{\footnotesize$\sin\frac{\varepsilon}{2}+\frac{x_1}{2-\cos\varepsilon}(\cos\frac{\varepsilon}{2}-\sin\frac{\varepsilon}{2})$}
\psfrag{-cose2}{\footnotesize$-\sin\frac{\varepsilon}{2}+\frac{x_1}{2-\cos\varepsilon}(\sin\frac{\varepsilon}{2}-\cos\frac{\varepsilon}{2})$}
\psfrag{sintt2}{\footnotesize$\sin\frac{\tilde{\theta}}{2}$}
\psfrag{Gs}{\footnotesize$\Gamma(\sigma)$}
\psfrag{-2}{\footnotesize$-2$}
\psfrag{0}{\footnotesize$0$}
\psfrag{2}{\footnotesize$2$}
\includegraphics[width=300pt]{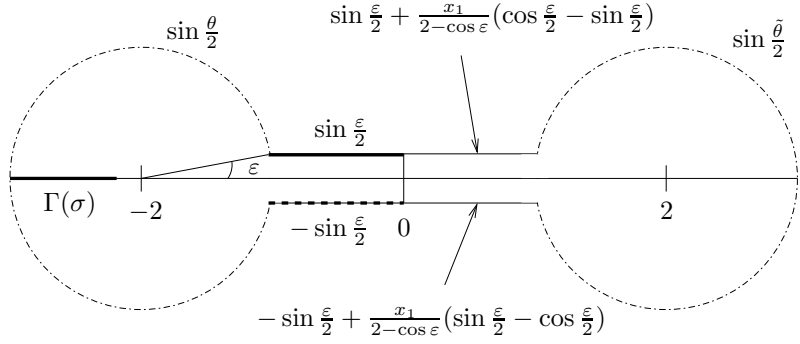}
\caption{\em The boundary datum $\psi_\varepsilon$.}
\end{center}
\end{figure}

For every $\sg \in ]-3,3[\,$, let $\ue(\sg)\in H^{1}(\Omega_\e\setminus \G(\sg))$ be
the  solution
of the problem:
\begin{equation}\llabel{eee}
E_\e(\sg):= \min\Big\{
\int_{\Omega_\e \setminus \G(\sg)} |D u|^2\, dx:u\in AD(\psi_\e,\sg)\Big\}.
\end{equation}
Our aim is to prove that for $\e$ sufficiently small there exists a subinterval 
$[a,b]$ of $[-2,2]$ such that
$E_\e(\sg)$ is concave on $[a,b]$.

As $\sg\mapsto E_\e(\sg)$ is a $C^2$-function, in order to prove that $E_\e(\sg)$
cannot be convex on the whole interval $[-2,2]$, it is enough to show that the
following three conditions are satisfied:
\begin{itemize}
\item[(a)] $\limsup_{\e\to 0^+} E_\e(2)$ is finite;
\item[(b)] $\liminf_{\e\to 0^+} E_\e(-2) = \infty$;
\item[(c)] $\limsup_{\e\to 0^+} E'_\e(-2)$ is finite;
\end{itemize}
where we denote by $'$ the first derivative with respect to $\sg$.

In order to prove condition (a) we construct an admissible function $\tilde u_\e$ for
$E_\e(2)$ whose energy, $\|D\tilde u_\e\|_2^2$, is bounded uniformly with respect to
$\e$. We define the open sets
$B^+_{-2}$ and
$B^-_{-2}$ by
\begin{equation*}
\begin{array}{c}
B^+_{-2}=\{(x_1,x_2)\in B_{-2}: x_2 > 0\}\\
B^-_{-2}=\{(x_1,x_2)\in B_{-2}: x_2 < 0\}.
\end{array}
\end{equation*}

Let $v^+$ 
be the solution to the following problem:
\begin{equation*}
\begin{cases}
\Delta u=0 &\mbox{on }B^+_{-2};\\
u(x)=\sin\frac{\theta(x)}{2} &\mbox{on }\partial B^+_{-2}\cap\partial B_{-2};\\
\partial_\nu u =0 &\mbox{on }{]-3,-1[}\times \{0\}.
\end{cases}
\end{equation*}
Then the function $v^-(x_1,x_2):= -v^+(x_1,-x_2)$ solves the analogue problem on
$B^-_{-2}$. 
Let $\tilde{u}_\e$ be the function which coincides with the harmonic functions
that satisfy the boundary conditions on $B^+_{-2}$, on $B^-_{-2}$, and on $B_2$,
respectively, that is, 
$\tilde{u}_\e:=v^+$ on $B^+_{-2}$,
$\tilde{u}_\e:=v^-$ on $B^-_{-2}$, and
$\tilde{u}_\e:= \tilde\rho^{\frac 1 2}\sin\frac{\tilde{\theta}}{2}$ on $B_2$. 
On $T_\e\setminus(B_2\cup B_{-2})$ we define $\tilde u_\e$ 
 in the following way: on the horizontal line $x_2=\sin\theta$, with
$\theta\in[-\e,\e]$, we set 
$\tilde{u}_\e(x_1,x_2):= \sin
\frac{\theta}{2}$ for $x_1\in ]-2+\cos\theta,0]$ and then interpolate linearly with
the boundary data on $\partial B_2\cap T_\e$: $\tilde{u}_\e(x_1,x_2):= 
\sin \frac{\theta}{2} +
\frac{x_1}{2-\cos\theta}(\cos\frac{\theta}{2}-\sin\frac{\theta}{2})$ for $x_1\in
[0,2-\cos\theta[\,$, if $0<\theta\leq\e$, and $\tilde{u}_\e(x_1,x_2):=  -\sin
\frac{\theta}{2} +
\frac{x_1}{2-\cos\theta}(-\sin\frac{\theta}{2}-\cos\frac{\theta}{2})$
for $x_1\in [0,2-\cos\theta[\,$, if $-\e\leq\theta<0$.
It is easy to check that $\tilde u_\e\in AD(\psi_\e,2)$ and that 
$D \tilde{u}_\e$ is bounded in $L^2({\Omega_\e\setminus \G};\R ^2)$ 
uniformly with respect to $\e$. 
This implies that
\begin{equation*}
\limsup_{\e\to 0^+} E_\e(2) \leq \limsup_{\e\to 0^+} 
\int_{\Omega_\e\setminus \G(2)}|D \tilde{u}_\e|^2\, dx < +\infty,
\end{equation*}
and condition (a) is satisfied.

We continue by proving condition (b), i.e., $E_\e(-2)$ tends to infinity as $\e$ goes
to zero. Let us first consider the model problem
\begin{equation}\llabel{ap}
\min\big\{\int_{R_\e}|D u|^2\, dx: u\geq \frac 1 2\, \mbox{ on }\partial_1 R_\e,
u\leq -\frac 1 2\, \mbox{ on }\partial_2R_\e\big\}
\end{equation}
where
\begin{equation}\llabel{Raux}
R_\e:= {]0,1[} \times {\big]-\e,\e\big[}, \quad
\partial_1R_\e := {[0,1]} \times \big\{\e \big\}, \quad
\partial_2R_\e := {[0,1]} \times \big\{- \e \big\}.
\end{equation}

 It is easy to see that problem (\ref{ap}) admits a
solution and that it is equivalent to
\begin{equation}\llabel{ap=}
\min\big\{\int_{R_\e}|D u|^2\, dx: u= \frac 1 2\, \mbox{ on }\partial_1R_\e, 
u= -\frac 1 2\, \mbox{ on }\partial_2R_\e\big\},
\end{equation}
which admits the affine solution $u^a(x_1,x_2):=\frac{1}{2\e}\, x_2$ 
for every $x=(x_1,x_2)\in R_\e$.

Going back to the domain $\Omega_\e$, let us consider the same problem with
different constants: the rectangle $R_\e$ is defined now by 
\begin{equation*}
R_\e:= {]A_\e,2-\cos\e[}\times{]-\sin\e ,\sin\e [}\subset T_\e\,,
\end{equation*}
where $A_\e$ is a positive constant such that 
$\psi_\e(x)\geq \frac 1 2$ on
$\partial_1  R_\e:=[A_\e,2-\cos\e]\times\{\sin\e \}$, (and
$\psi_\e(x)\leq -\frac 1 2$ on
$\partial_2  R_\e:=[A_\e,2-\cos\e]\times\{-\sin\e \}$), when $\e$ is
sufficiently small. 
Then 
\begin{equation*}
E_\e(-2)=\int_{\Om_\e\setminus\Gamma(-2)}|Du^\e(-2)|^2dx\geq
\int_{R_\e}|Du^\e(-2)|^2dx\geq\int_{R_\e}|D u^a|^2 dx\,.
\end{equation*}  
Since  $\int_{R_\e}|D u^a|^2 dx\to+\infty$ as $\e\to0$,
condition (b) is proved.

It remains to show that condition (c) is satisfied, i.e., that the first
derivative of 
$\sg\mapsto E_\e(\sg)$ at $\sg=-2$ is bounded  as $\e$ goes to zero.
Since 
\begin{equation}\llabel{dea}
E'_\e(\sg) = - \kk_\e^2(\sg) \,,
\end{equation}
see, e.g. \cite[Theorem 6.4.1]{Gr-92}, where $\kk_\e(\sg)\sqrt{\frac{2}{\pi}}$ is
the stress intensity factor associated to $u^\e(\sg)$ at the tip $(\sg,0)$,
see Proposition~\ref{th:2.1}, it is enough to show that $\kk_\e(\sg)$ remains
bounded when, for instance, $-\frac 5 2\leq\sg\leq -\frac 3 2$.

For $\sg\in [-5/2,-3/2]$, let 
$v(\sg)$ be the solution of the following problem:
\begin{equation}\llabel{v-2}
\min\big\{ \int_{B_{-2}\setminus \G(\sg)}|D u|^2\, dx: 
u\in H^{1}(B_{-2}\setminus \G(\sg)),\,u=\sin\frac{\theta}{2}\,\mbox{ on }\partial
B_{-2}\setminus \G(\sg)\big\}.
\end{equation}
Let us extend $v(\sg)$ to $\R \times[-1,1]$ constantly on the horizontal lines and 
denote now by $v(\sg)$ this extension.% (which out of $B_{-2}$ depends only on~$x_2$).

We claim that
\begin{equation}\llabel{claim}
u^\e(\sg)\to v(\sg) \quad \mbox{strongly in }H^{1}(B_{-2}\setminus \G(\sg)).
\end{equation}

Assuming the claim true, we now use the following characterization of
$\kk_\e$ (see Proposition \ref{p:kchar}):
\begin{equation}\llabel{ka}
\kk_\e^2(\sg) = \int_{B_{-2}\setminus\G(\sg)} 
\Big[ ((D_1 u^\e)^2-(D_2 u^\e)^2)D_1 \varphi + 2D_1u^\e D_2 u^\e D_2\varphi\Big]\, dx 
\end{equation}
with $\varphi\in C^1_c(B_{-2})$ such that 
$\varphi(\sg,0) = 1$. 
By (\ref{claim}) and the definition of $v(\sg)$, we can pass to the limit in the right-hand side as 
$\e\to 0^+$ and define in such a way the quantity:
\begin{equation}\llabel{ksg}
\kk^2(\sg) := \int_{B_{-2}\setminus\G(\sg)} 
\Big[ ((D_1 v(\sg))^2-(D_2 v(\sg))^2)D_1 \varphi + 2D_1v(\sg)D_2 v(\sg)D_2\varphi\Big]\, dx.
\end{equation}
Therefore, by (\ref{dea}), 
\begin{equation}\llabel{c0}
\limsup_{\e\to 0^+} E'_\e(\sg) = -\kk^2(\sg) \quad\mbox{for every }-
\frac 5 2\leq\sg\leq -\frac 3 2\,.
\end{equation}
As, by (\ref{ksg}), $\kk(\sg)$ is bounded, formula (\ref{c0}) concludes the 
proof of condition (c).

\noindent{\em Proof of the claim. }
Let $\tilde{\Omega}_\e:=T_\e \cup B_2$ 
and let $w_\e$ be the solution of the following problem:
\begin{equation}\llabel{wA}
\min\Big\{ \int_{\tilde{\Omega}_\e}|D u|^2\, dx: 
u\in H^{1}(\tilde{\Omega}_\e),\,u=\psi_\e\,\mbox{ on }\partial\Omega_\e
\cap\partial\tilde{\Omega}_\e\Big\}.
\end{equation}
We consider a cut-off function $\varphi\in C^\infty(\R )$ such that 
$0\leq \varphi\leq 1$, $\varphi(x_1)=1$ for $x_1\leq -\frac 2 3$, and 
$\varphi(x_1)=0$ for $x_1\geq -\frac 1 3$.
Then the function $\zeta:=\varphi v(\sg) + (1-\varphi)w_\e$ belongs to
$AD(\psi_\e,\sg)$ 
and 
\begin{equation}\llabel{mini}
E_\e(\sg)=\int_{\Omega_\e\setminus \G(\sg)}|D u^\e(\sg)|^2\, dx
\leq \int_{\Omega_\e\setminus \G(\sg)}|D \zeta|^2\, dx.
\end{equation}
By convexity, we have 
\begin{equation}\llabel{keyaux}
\begin{array}{c}
{\displaystyle \int_{\Omega_\e\setminus \G(\sg)}|D \zeta|^2\, dx\leq
\int_{B_{-2}\setminus\G(\sg)}|D v(\sg)|^2\, dx + 
\int_{\tilde{\Omega}_\e}|D w_\e|^2\, dx +\int_{T_\e} |D v(\sg)|^2\, dx +}\\[3mm]
{\displaystyle + \int_{T_\e\cap ({\rm supp} D\varphi)} \big(2D\varphi
(\varphi D v(\sg) + (1-\varphi)D w_\e)(v(\sg)-w_\e)+
|D\varphi|^2(v(\sg)-w_\e)^2\big)\, dx}.
\end{array}
\end{equation}
Now 
\begin{equation}\llabel{i}
\begin{array}{l}
{\displaystyle \lim_{\e\to 0^+}\int_{T_\e}|D v(\sg)|^2\, dx = 0,}\\[3mm]
{\displaystyle \lim_{\e\to 0^+}\int_{T_\e\cap ({\rm supp} D\varphi)}
|D\varphi|^2(v(\sg)-w_\e)^2\, dx = 0,}
\end{array}
\end{equation}
and, for any $\eta>0$,
\begin{equation}\llabel{ii}
\begin{array}{c} 
{\displaystyle \int_{T_\e\cap ({\rm supp} D\varphi)} 
2D\varphi(\varphi D v(\sg) + (1-\varphi)D w_\e)(v(\sg)-w_\e)\, dx
\leq }\\
{\displaystyle \leq 2\int_{T_\e\cap ({\rm supp} D\varphi)} 
D\varphi \varphi D v(\sg) (v(\sg)-w_\e)\, dx +}\\
{\displaystyle + \frac{1}{\eta}\int_{T_\e\cap ({\rm supp}
D\varphi)}|D\varphi|^2|v(\sg)-w_\e|^2\, dx+ \eta
\int_{T_\e\cap ({\rm supp} D\varphi)}
|D w_\e|^2(1-\varphi)^2\, dx}\,.
\end{array}
\end{equation}
Since the first two terms in the right-hand side tend to zero, it remains to prove
that 
\begin{equation}\llabel{wa} 
\lim_{\e\to 0^+}\int_{T_\e\cap ({\rm supp} D\varphi)}|D w_\e|^2\, dx=0\,.
\end{equation}
As in the proof of condition (b), we  consider first a model problem. 
Similarly to (\ref{Raux}), we now set
$$
R_\e:=]-1,0[\,\times\,]-\e ,\e[\,,\quad
\partial_1R_\e:=[-1,0]\times\{\e \},\quad
\partial_2R_\e:=[-1,0]\times\{-\e \},
$$
and define $h_\e$ as the solution to the following problem:
\begin{equation}\llabel{ha}
\begin{cases}
\Delta h_\e = 0 &\mbox{ on }R_\e,\\
h_\e=\frac \e 2  &\mbox{ on }\partial_1 R_\e,\\
h_\e=-\frac \e 2 &\mbox{ on }\partial_2 R_\e,\\
\|h_\e\|_\infty \leq 1.
\end{cases}
\end{equation}
We claim that
\begin{equation}\llabel{claim2}
\lim_{\e\to 0^+} \int_{\tilde{R}_\e} |D h_\e|^2\, dx = 0,
\end{equation}
where
\begin{equation*}
\tilde{R}_\e := {\big]-\frac 4 5,-\frac 1 5\big[}\times {]- \e,\e [}\subset R_\e.
\end{equation*}
Indeed, note that the function $z_\e(x_1,x_2):=\frac 1 2 x_2$ solves (\ref{ha}) (for
$\e\leq1$).
By a Cacciopoli type estimate we obtain
$$
\int_{\tilde R_\e}|D(h_\e-z_\e)|^2dx\leq C\int_{R_\e}|h_\e-z_\e|^2dx\leq C_1|R_\e|\,,
$$
for some positive constants $C$ and $C_1$ which do not depend on $\e$, hence
(\ref{claim2}) holds.

Applying this argument with
\begin{equation*}
R_\e={\big]-1,0\big[}\times {]-\sin\e ,\sin\e [}\quad\mbox{ and }
\quad\tilde{R}_\e={\big]-\frac 4 5,-\frac 1 5\big[}\times {]-\sin\e,\sin\e[}
\end{equation*}
it follows that (\ref{wa}) holds true.
\begin{figure}[ht]
\begin{center}
\psfrag{B-}{$B_{-2}$}
\psfrag{Gs}{$\Gamma(\sigma)$}
\psfrag{-2}{$-2$}
\psfrag{0}{$0$}
\psfrag{tRe}{$\tilde{R}_\epsilon$}
\psfrag{tOe}{$\tilde{\Omega}_\epsilon$}
\psfrag{2}{$2$}
\psfrag{}{$$}
\includegraphics[width=300pt]{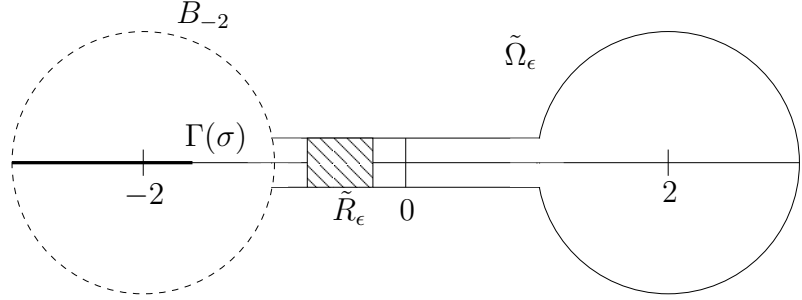}
\caption{\em The rectangle $\tilde{R}_\epsilon$ where we apply a Cacciopoli type estimate in order to obtain (\ref{wa}).}
\end{center}
\end{figure}

{}From (\ref{mini}), (\ref{keyaux}), (\ref{i}), (\ref{ii}), %(\ref{iii}), 
and (\ref{wa}) we deduce that 
\begin{equation}\llabel{fond}
\int_{\Omega_\e\setminus \G(\sg)}|D u^\e(\sg)|^2\,dx \leq \int_{B_{-2}\setminus
\G(\sg)}|D v(\sg)|^2\, dx +
\int_{\tilde{\Omega}_\e}|D w_\e|^2\, dx + o(1).
\end{equation}
Since 
$$
\int_{\tilde{\Omega}_\e}|Du^\e(\sg)|^2dx \geq\int_{\tilde{\Omega}_\e}|D w_\e|^2 dx 
$$
we obtain 
\begin{equation}\llabel{fond1}
\int_{B_{-2}\setminus \G(\sg)}|D u_\e(\sg)|^2\,dx \leq 
\int_{B_{-2}\setminus \G(\sg)}|D v(\sg)|^2\, dx + o(1)\leq C
\end{equation} 
 uniformly with respect
to $\e$. Thus, there exists $u^*(\sg)\in H^{1}(B_{-2}\setminus \G(\sg))$ such that 
\begin{equation}\llabel{u*}
u^\e(\sg)\rightharpoonup u^*(\sg)\quad \mbox{weakly on }H^{1}(B_{-2}\setminus
\G(\sg)),
\end{equation}
and
\begin{equation}\llabel{bc*}
u^*(\sg)=\sin\frac{\theta}{2} \quad \mbox{on }\partial B_{-2}\setminus \G(\sg).
\end{equation}
As $(D u^\e(\sg),D\varphi)=0$ for every $\varphi\in H^1(B_{-2}\setminus \G(\sg))$
with $\varphi=0 $ on $\partial B_{-2}\setminus \G(\sg)$, by (\ref{u*}) we obtain
that $(D u^*(\sg),D \varphi)=0$. By (\ref{v-2}), this fact, together with
(\ref{bc*}), implies that 
\begin{equation}\llabel{u=v}
u^*(\sg)=v(\sg)\,.
\end{equation}  
In addition, by the lower semicontinuity and by (\ref{fond1}), we have 
\begin{equation}\llabel{strong}
\int_{B_{-2}\setminus \G(\sg)}|D v(\sg)|^2\, dx \leq \liminf_{\e\to 0^+}
\int_{B_{-2}\setminus \G(\sg)}|D u^\e(\sg)|^2\, dx \leq \int_{B_{-2}\setminus
\G(\sg)}|D v(\sg)|^2\, dx.
\end{equation}
By (\ref{u*}), (\ref{u=v}), and (\ref{strong}), we deduce that (\ref{claim}) holds.

\end{section}

\bigskip
\bigskip

\noindent
{{\bf Acknowledgements.}
\nopagebreak  
\medskip

\noindent
The authors wish to thank Gianni Dal Maso for many interesting and fruitful
discussions. This work is part of the project ``Calculus of Variations" 2004  supported by the Italian Ministry of Education, University, and Research.}

\end{document}